\documentclass{amsproc}
\usepackage{color}

\newtheorem{theorem}{Theorem}[section]
\newtheorem{lemma}[theorem]{Lemma}

\newtheorem{corollary}[theorem]{Corollary}
\theoremstyle{definition}
\newtheorem{definition}[theorem]{Definition}

\theoremstyle{remark}
\newtheorem{remark}[theorem]{Remark}

\newtheorem{notation}[theorem]{Notation}

\numberwithin{equation}{section}

\newcommand{\abs}[1]{\lvert#1\rvert}
\newcommand{\Abs}[1]{\Vert#1\Vert}
\newcommand{\ABS}[1]{\Vert\lvert#1\rvert\Vert}
\newcommand{\blankbox}[2]{%
  \parbox{\columnwidth}{\centering

    \setlength{\fboxsep}{0pt}%
    \fbox{\raisebox{0pt}[#2]{\hspace{#1}}}%
  }%
}

\newcommand{\s}{{\rm supp}}
\newcommand{\R}{{\mathbb R}}
\newcommand{\RR}{\mathcal R}
\newcommand{\Z}{{\mathbb Z}}
\newcommand{\FF}{{\mathcal F}}
\newcommand{\GG}{{\mathcal G}}
\newcommand{\Int}{{{\rm Int}\,}}
\newcommand{\A}{{\mathbb A}}
\newcommand{\N}{{\mathbb N}}
\newcommand{\LL}{{\mathcal L}}
\newcommand{\IA}{{\stackrel{\circ}{\A}}}
\newcommand{\C}{{\mathbb C}}
\newcommand{\Fix}{{\rm Fix}}
\newcommand{\Moeb}{\mbox{M\"ob}_+(S^1_\C)}
\newcommand{\BB}{{\mathcal B}}
\newcommand{\PP}{{\mathbb P}}
\newcommand{\EE}{{\mathbb E}}
\newcommand{\HH}{{\mathcal H}}
\newcommand{\lev}{{\rm lev}}
\newcommand{\Stab}{{\rm Stab}}
\newcommand{\bd}{{\sc Proof}.\ }
\newcommand{\QQ}{{\mathcal Q}}
\newcommand{\VV}{{\mathbb V}}
\newcommand{\JJ}{{\mathcal J}}
\newcommand{\TT}{{\mathcal T}}
\newcommand{\II}{{\stackrel{\circ} I}}
\newcommand{\VVV}{{\mathcal V}}
\newcommand{\EEE}{{\mathcal E}}
\newcommand{\DD}{{\mathbb D}}

\begin{document}
\title[Planar Anosov diffeomorphisms]{An example of planar Anosov diffeomorphisms without fixed points}


\author{Shigenori Matsumoto}
\address{Department of Mathematics, College of
Science and Technology, Nihon University, 1-8-14 Kanda, Surugadai,
Chiyoda-ku, Tokyo, 101-8308 Japan
}
\email{matsumo@math.cst.nihon-u.ac.jp
}
\thanks{The author is partially supported by Grant-in-Aid for
Scientific Research (C) No.\ 18K03312.}


\date{\today }
\begin{abstract}
We construct an example of fixed point free
Anosov diffeomorphisms of the plane, which is not topological
conjugate to a translation.
\end{abstract}

\maketitle
\section{Introduction}
In \cite{W}, W. White showed that a translation of the plane is an
Anosov diffeomorphism in the sense of Definition \ref{def}.
P. Mendes \cite{M} studied properties of Anosov diffeomorphisms
 of the plane and 
 conjectured that any planar fixed point free Anosov
diffeomorphism is topologically conjugate to a translation. 
The purpose of this paper is to disprove this conjecture. 
First let us recall the definition of Anosov diffeomorphisms of
the plane.

\begin{definition}\label{def}
A $C^1$ diffeomorphism $F$ of the plane $\R^2$ is said to be an {\em Anosov
 diffeomorphism} if there are a continuous Riemannian metric 
$m$ and two transversal continuous foliations $\FF^u$ and $\FF^s$
by $C^1$-leaves with the following
properties:
\\
(1) the metric $m$ is complete,
\\
(2) the diffeomorphism $F$ preserves the two foliations $\FF^\sigma$, $\sigma=u,s$,
i.e, maps each leaf of $\FF^\sigma$ to a leaf of $\FF^\sigma$, and
\\
(3) there are constant $C>0$ and $\lambda>0$ such that
\begin{equation}\label{1}
\Abs{DF^n(v)}_m\geq C^{-1}e^{\lambda n}\Abs{v}_m, \ \forall v\in T\FF^u,
\ \forall n\in \N,
\end{equation}
and
\begin{equation}\label{2}
\Abs{DF^n(v)}_m\leq Ce^{-\lambda n}\Abs{v}_m, \forall
v\in T\FF^s\ \forall n\in\N.
\end{equation}
\end{definition}
\smallskip
The condition (1) is necessary in order to exclude trivial
examples. Consider a linear diffeomorphism
$A$ defined by $A(x,y)=(2x,\tfrac{1}{2}y)$ and consider an $A$-invariant
strip 
$$C=\{(x,y)\in\R^2\mid x>0, 1<xy<2\}.$$
Then $A\vert_C:C\to C$ satisfies conditions (2) and (3) with respect
to the vertical and horizontal foliations, and the
 metric $m$ which is the restriction of the Euclidean metric to
$C$. However $m$ is not complete. The example in \cite{W} is more involved.

\medskip
{\sc Main Theorem}. \em
There is 
a fixed point free Anosov diffeomorphism which is not topologically
conjugate to a translation. 
\rm 

\medskip

J. Groisman and Z. Nitecki \cite{GN} proved the Mendes conjecture for
a certain class of diffeomorphisms i.e.\ the time one maps of
 $C^1$-flows.
In fact, they showed the following.

\begin{theorem}\label{t1}
Let $F$ be the time one map of a fixed point free $C^1$ flow which is not
 topological conjugate to a translation. Assume $F$ preserves a continuous
 foliation $\FF$ by $C^1$ leaves. Then some leaf $L$ of $\FF$ is left invariant by $F$.
\end{theorem}
This quickly leads to the solution of the Mendes conjecture for this
class
of diffeomorphisms, since if 
$\FF=\FF^u$, $F$ must have a fixed point
in $L$ by virtue of (\ref{1}).

Therefore our first task for the proof of Main Theorem
is to construct a $C^1$ diffeomorphism $F$
and two mutually transverse foliations, say $\FF^u$ and $\FF^s$,
 invariant by $F$ but without invariant leaves. The schematic idea
can be found in Figure 1.
The solid lines indicate the foliation $\FF^u$, while dotted lines
$\FF^s$. The diffeomorphism $F$ maps $p_i$ to $p_{i+1}$, and
$q_i$ to $q_{i-1}$.
 Detailed construction is described in Sections 2 and 3.
It may be worth mentioning that there is no contradiction with the
Brouwer plane fixed point theorem. Horizontal and vertical ``Reeb
components'' are displaced, and outside them, the diffeomorphism is
conjugate to a translation of the plane. Thus all the points are
wandering.
Sections 4 and 5 are devoted to
  the definition of the metric.

\begin{figure}[h]
{\unitlength 0.1in%
\begin{picture}( 49.4600, 27.5800)( 14.9000,-87.8100)%
%
\special{pn 4}%
\special{sh 1}%
\special{ar 1735 7017 16 16 0  6.28318530717959E+0000}%
\special{sh 1}%
\special{ar 1735 7017 16 16 0  6.28318530717959E+0000}%
%
\special{pn 4}%
\special{sh 1}%
\special{ar 2527 6232 16 16 0  6.28318530717959E+0000}%
\special{sh 1}%
\special{ar 2527 6232 16 16 0  6.28318530717959E+0000}%
%
\special{pn 4}%
\special{sh 1}%
\special{ar 2282 7593 16 16 0  6.28318530717959E+0000}%
\special{sh 1}%
\special{ar 2282 7593 16 16 0  6.28318530717959E+0000}%
%
\special{pn 4}%
\special{sh 1}%
\special{ar 2851 8162 16 16 0  6.28318530717959E+0000}%
\special{sh 1}%
\special{ar 2851 8162 16 16 0  6.28318530717959E+0000}%
%
\special{pn 4}%
\special{sh 1}%
\special{ar 3074 6786 16 16 0  6.28318530717959E+0000}%
\special{sh 1}%
\special{ar 3074 6786 16 16 0  6.28318530717959E+0000}%
%
\special{pn 4}%
\special{sh 1}%
\special{ar 3664 7370 16 16 0  6.28318530717959E+0000}%
\special{sh 1}%
\special{ar 3672 7362 16 16 0  6.28318530717959E+0000}%
%
\special{pn 4}%
\special{sh 1}%
\special{ar 4240 7946 16 16 0  6.28318530717959E+0000}%
\special{sh 1}%
\special{ar 4240 7938 16 16 0  6.28318530717959E+0000}%
%
\special{pn 8}%
\special{pa 3002 6081}%
\special{pa 4788 6081}%
\special{fp}%
%
\special{pn 8}%
\special{pa 3002 6354}%
\special{pa 4788 6354}%
\special{fp}%
%
\special{pn 8}%
\special{pa 3002 6088}%
\special{pa 2972 6087}%
\special{pa 2941 6086}%
\special{pa 2908 6087}%
\special{pa 2872 6089}%
\special{pa 2834 6094}%
\special{pa 2794 6100}%
\special{pa 2754 6109}%
\special{pa 2714 6119}%
\special{pa 2675 6130}%
\special{pa 2639 6142}%
\special{pa 2606 6156}%
\special{pa 2577 6170}%
\special{pa 2554 6185}%
\special{pa 2538 6200}%
\special{pa 2529 6215}%
\special{pa 2528 6231}%
\special{pa 2536 6246}%
\special{pa 2551 6261}%
\special{pa 2573 6275}%
\special{pa 2601 6289}%
\special{pa 2634 6302}%
\special{pa 2670 6314}%
\special{pa 2708 6325}%
\special{pa 2748 6334}%
\special{pa 2788 6342}%
\special{pa 2828 6348}%
\special{pa 2866 6353}%
\special{pa 2902 6355}%
\special{pa 2935 6356}%
\special{pa 2966 6355}%
\special{pa 2996 6354}%
\special{pa 3002 6354}%
\special{fp}%
%
\special{pn 8}%
\special{pa 3556 6657}%
\special{pa 3526 6656}%
\special{pa 3495 6655}%
\special{pa 3462 6656}%
\special{pa 3426 6658}%
\special{pa 3388 6663}%
\special{pa 3348 6669}%
\special{pa 3308 6678}%
\special{pa 3268 6688}%
\special{pa 3229 6699}%
\special{pa 3193 6711}%
\special{pa 3160 6725}%
\special{pa 3131 6739}%
\special{pa 3108 6754}%
\special{pa 3092 6769}%
\special{pa 3083 6784}%
\special{pa 3082 6800}%
\special{pa 3090 6815}%
\special{pa 3105 6830}%
\special{pa 3127 6844}%
\special{pa 3155 6858}%
\special{pa 3188 6871}%
\special{pa 3224 6883}%
\special{pa 3262 6894}%
\special{pa 3302 6903}%
\special{pa 3342 6911}%
\special{pa 3382 6917}%
\special{pa 3420 6922}%
\special{pa 3456 6924}%
\special{pa 3489 6925}%
\special{pa 3520 6924}%
\special{pa 3550 6923}%
\special{pa 3556 6923}%
\special{fp}%
%
\special{pn 8}%
\special{pa 4125 7240}%
\special{pa 4095 7239}%
\special{pa 4064 7238}%
\special{pa 4031 7239}%
\special{pa 3995 7241}%
\special{pa 3957 7246}%
\special{pa 3917 7252}%
\special{pa 3877 7261}%
\special{pa 3837 7271}%
\special{pa 3798 7282}%
\special{pa 3762 7294}%
\special{pa 3729 7308}%
\special{pa 3700 7322}%
\special{pa 3677 7337}%
\special{pa 3661 7352}%
\special{pa 3652 7367}%
\special{pa 3651 7383}%
\special{pa 3659 7398}%
\special{pa 3674 7413}%
\special{pa 3696 7427}%
\special{pa 3724 7441}%
\special{pa 3757 7454}%
\special{pa 3793 7466}%
\special{pa 3831 7477}%
\special{pa 3871 7486}%
\special{pa 3911 7494}%
\special{pa 3951 7500}%
\special{pa 3989 7505}%
\special{pa 4025 7507}%
\special{pa 4058 7508}%
\special{pa 4089 7507}%
\special{pa 4119 7506}%
\special{pa 4125 7506}%
\special{fp}%
%
\special{pn 8}%
\special{pa 4701 7802}%
\special{pa 4671 7801}%
\special{pa 4640 7800}%
\special{pa 4607 7801}%
\special{pa 4571 7803}%
\special{pa 4533 7808}%
\special{pa 4493 7814}%
\special{pa 4453 7822}%
\special{pa 4413 7832}%
\special{pa 4374 7843}%
\special{pa 4338 7856}%
\special{pa 4305 7869}%
\special{pa 4276 7883}%
\special{pa 4253 7898}%
\special{pa 4237 7913}%
\special{pa 4227 7928}%
\special{pa 4227 7944}%
\special{pa 4235 7959}%
\special{pa 4250 7974}%
\special{pa 4273 7989}%
\special{pa 4300 8003}%
\special{pa 4333 8016}%
\special{pa 4369 8028}%
\special{pa 4407 8039}%
\special{pa 4447 8048}%
\special{pa 4487 8056}%
\special{pa 4527 8062}%
\special{pa 4565 8067}%
\special{pa 4601 8069}%
\special{pa 4634 8070}%
\special{pa 4665 8069}%
\special{pa 4695 8068}%
\special{pa 4701 8068}%
\special{fp}%
%
\special{pn 8}%
\special{pa 3549 6657}%
\special{pa 5335 6657}%
\special{fp}%
%
\special{pn 8}%
\special{pa 3549 6916}%
\special{pa 5335 6916}%
\special{fp}%
%
\special{pn 8}%
\special{pa 4104 7247}%
\special{pa 5889 7247}%
\special{fp}%
%
\special{pn 8}%
\special{pa 4068 7506}%
\special{pa 5853 7506}%
\special{fp}%
%
\special{pn 8}%
\special{pa 4651 7802}%
\special{pa 6436 7802}%
\special{fp}%
%
\special{pn 8}%
\special{pa 4644 8068}%
\special{pa 6429 8068}%
\special{fp}%
%
\special{pn 8}%
\special{pa 1612 7449}%
\special{pa 1612 8781}%
\special{dt 0.045}%
%
\special{pn 8}%
\special{pa 1836 7442}%
\special{pa 1850 8752}%
\special{dt 0.045}%
%
\special{pn 8}%
\special{pn 8}%
\special{pa 1612 7449}%
\special{pa 1611 7441}%
\special{fp}%
\special{pa 1607 7403}%
\special{pa 1606 7394}%
\special{fp}%
\special{pa 1604 7356}%
\special{pa 1604 7348}%
\special{fp}%
\special{pa 1606 7310}%
\special{pa 1607 7301}%
\special{fp}%
\special{pa 1612 7263}%
\special{pa 1613 7255}%
\special{fp}%
\special{pa 1620 7218}%
\special{pa 1622 7209}%
\special{fp}%
\special{pa 1631 7172}%
\special{pa 1633 7164}%
\special{fp}%
\special{pa 1645 7128}%
\special{pa 1648 7120}%
\special{fp}%
\special{pa 1664 7085}%
\special{pa 1668 7078}%
\special{fp}%
\special{pa 1690 7047}%
\special{pa 1696 7041}%
\special{fp}%
\special{pa 1731 7035}%
\special{pa 1739 7038}%
\special{fp}%
\special{pa 1764 7067}%
\special{pa 1769 7073}%
\special{fp}%
\special{pa 1787 7107}%
\special{pa 1790 7115}%
\special{fp}%
\special{pa 1804 7150}%
\special{pa 1807 7158}%
\special{fp}%
\special{pa 1818 7195}%
\special{pa 1819 7203}%
\special{fp}%
\special{pa 1827 7240}%
\special{pa 1829 7249}%
\special{fp}%
\special{pa 1834 7287}%
\special{pa 1835 7295}%
\special{fp}%
\special{pa 1837 7333}%
\special{pa 1837 7341}%
\special{fp}%
\special{pa 1835 7380}%
\special{pa 1834 7388}%
\special{fp}%
\special{pa 1829 7426}%
\special{pa 1828 7434}%
\special{fp}%
%
\special{pn 8}%
\special{pn 8}%
\special{pa 2174 8003}%
\special{pa 2173 7995}%
\special{fp}%
\special{pa 2169 7957}%
\special{pa 2168 7949}%
\special{fp}%
\special{pa 2166 7910}%
\special{pa 2166 7902}%
\special{fp}%
\special{pa 2168 7864}%
\special{pa 2169 7856}%
\special{fp}%
\special{pa 2174 7818}%
\special{pa 2175 7810}%
\special{fp}%
\special{pa 2182 7772}%
\special{pa 2184 7764}%
\special{fp}%
\special{pa 2193 7727}%
\special{pa 2195 7719}%
\special{fp}%
\special{pa 2208 7683}%
\special{pa 2211 7675}%
\special{fp}%
\special{pa 2226 7640}%
\special{pa 2230 7633}%
\special{fp}%
\special{pa 2252 7601}%
\special{pa 2258 7596}%
\special{fp}%
\special{pa 2293 7590}%
\special{pa 2300 7594}%
\special{fp}%
\special{pa 2326 7622}%
\special{pa 2331 7629}%
\special{fp}%
\special{pa 2348 7663}%
\special{pa 2352 7670}%
\special{fp}%
\special{pa 2365 7706}%
\special{pa 2367 7714}%
\special{fp}%
\special{pa 2379 7750}%
\special{pa 2381 7759}%
\special{fp}%
\special{pa 2389 7796}%
\special{pa 2391 7804}%
\special{fp}%
\special{pa 2396 7842}%
\special{pa 2397 7850}%
\special{fp}%
\special{pa 2399 7888}%
\special{pa 2399 7896}%
\special{fp}%
\special{pa 2396 7935}%
\special{pa 2395 7943}%
\special{fp}%
\special{pa 2391 7981}%
\special{pa 2390 7989}%
\special{fp}%
%
\special{pn 8}%
\special{pn 8}%
\special{pa 2743 8601}%
\special{pa 2742 8593}%
\special{fp}%
\special{pa 2738 8555}%
\special{pa 2737 8546}%
\special{fp}%
\special{pa 2735 8508}%
\special{pa 2735 8500}%
\special{fp}%
\special{pa 2737 8462}%
\special{pa 2738 8454}%
\special{fp}%
\special{pa 2743 8416}%
\special{pa 2744 8407}%
\special{fp}%
\special{pa 2751 8370}%
\special{pa 2753 8362}%
\special{fp}%
\special{pa 2762 8325}%
\special{pa 2764 8317}%
\special{fp}%
\special{pa 2776 8280}%
\special{pa 2779 8272}%
\special{fp}%
\special{pa 2795 8238}%
\special{pa 2799 8230}%
\special{fp}%
\special{pa 2820 8199}%
\special{pa 2827 8193}%
\special{fp}%
\special{pa 2861 8187}%
\special{pa 2869 8191}%
\special{fp}%
\special{pa 2894 8219}%
\special{pa 2899 8225}%
\special{fp}%
\special{pa 2917 8259}%
\special{pa 2920 8267}%
\special{fp}%
\special{pa 2934 8302}%
\special{pa 2937 8310}%
\special{fp}%
\special{pa 2948 8347}%
\special{pa 2949 8355}%
\special{fp}%
\special{pa 2957 8393}%
\special{pa 2959 8401}%
\special{fp}%
\special{pa 2964 8439}%
\special{pa 2965 8447}%
\special{fp}%
\special{pa 2967 8485}%
\special{pa 2967 8493}%
\special{fp}%
\special{pa 2965 8532}%
\special{pa 2964 8540}%
\special{fp}%
\special{pa 2960 8578}%
\special{pa 2959 8586}%
\special{fp}%
%
\special{pn 8}%
\special{pa 2174 7982}%
\special{pa 2174 8752}%
\special{dt 0.045}%
%
\special{pn 8}%
\special{pa 2390 7996}%
\special{pa 2390 8759}%
\special{dt 0.045}%
%
\special{pn 8}%
\special{pa 2728 8586}%
\special{pa 2743 8738}%
\special{dt 0.045}%
%
\special{pn 8}%
\special{pa 2959 8586}%
\special{pa 2966 8759}%
\special{dt 0.045}%
%
\special{pn 8}%
\special{pa 1735 6405}%
\special{pa 1735 8766}%
\special{fp}%
%
\special{pn 8}%
\special{pa 2282 7204}%
\special{pa 2282 8745}%
\special{fp}%
%
\special{pn 8}%
\special{pa 2851 7766}%
\special{pa 2844 8759}%
\special{fp}%
\special{pa 2844 8759}%
\special{pa 2844 8759}%
\special{fp}%
%
\special{pn 8}%
\special{pa 3441 8716}%
\special{pa 3448 8241}%
\special{fp}%
%
\special{pn 8}%
\special{pa 3448 8270}%
\special{pa 3442 8237}%
\special{pa 3438 8204}%
\special{pa 3435 8172}%
\special{pa 3436 8140}%
\special{pa 3440 8109}%
\special{pa 3449 8078}%
\special{pa 3464 8049}%
\special{pa 3483 8020}%
\special{pa 3506 7992}%
\special{pa 3530 7964}%
\special{pa 3556 7938}%
\special{pa 3581 7912}%
\special{pa 3605 7886}%
\special{pa 3626 7862}%
\special{pa 3644 7838}%
\special{pa 3657 7814}%
\special{pa 3663 7791}%
\special{pa 3663 7768}%
\special{pa 3655 7746}%
\special{pa 3641 7724}%
\special{pa 3622 7701}%
\special{pa 3599 7679}%
\special{pa 3573 7657}%
\special{pa 3545 7633}%
\special{pa 3516 7610}%
\special{pa 3487 7585}%
\special{pa 3460 7560}%
\special{pa 3435 7534}%
\special{pa 3414 7506}%
\special{pa 3398 7477}%
\special{pa 3387 7446}%
\special{pa 3382 7414}%
\special{pa 3384 7380}%
\special{pa 3390 7346}%
\special{pa 3400 7311}%
\special{pa 3412 7277}%
\special{pa 3423 7243}%
\special{pa 3434 7211}%
\special{pa 3441 7181}%
\special{pa 3445 7154}%
\special{pa 3442 7130}%
\special{pa 3433 7109}%
\special{pa 3417 7092}%
\special{pa 3394 7078}%
\special{pa 3367 7066}%
\special{pa 3336 7055}%
\special{pa 3303 7046}%
\special{pa 3267 7037}%
\special{pa 3231 7027}%
\special{pa 3194 7017}%
\special{pa 3159 7005}%
\special{pa 3125 6991}%
\special{pa 3095 6975}%
\special{pa 3069 6955}%
\special{pa 3047 6931}%
\special{pa 3027 6906}%
\special{pa 3024 6902}%
\special{fp}%
%
\special{pn 8}%
\special{pa 4010 8730}%
\special{pa 3994 8666}%
\special{pa 3987 8634}%
\special{pa 3982 8602}%
\special{pa 3979 8570}%
\special{pa 3978 8538}%
\special{pa 3981 8506}%
\special{pa 3987 8474}%
\special{pa 3997 8442}%
\special{pa 4012 8410}%
\special{pa 4032 8378}%
\special{pa 4054 8347}%
\special{pa 4077 8316}%
\special{pa 4098 8287}%
\special{pa 4115 8258}%
\special{pa 4126 8232}%
\special{pa 4129 8207}%
\special{pa 4122 8184}%
\special{pa 4104 8164}%
\special{pa 4078 8145}%
\special{pa 4048 8126}%
\special{pa 4017 8107}%
\special{pa 3989 8087}%
\special{pa 3967 8064}%
\special{pa 3954 8039}%
\special{pa 3951 8011}%
\special{pa 3956 7980}%
\special{pa 3968 7947}%
\special{pa 3984 7914}%
\special{pa 4004 7880}%
\special{pa 4025 7847}%
\special{pa 4045 7815}%
\special{pa 4064 7784}%
\special{pa 4078 7756}%
\special{pa 4087 7731}%
\special{pa 4089 7710}%
\special{pa 4082 7693}%
\special{pa 4068 7680}%
\special{pa 4047 7669}%
\special{pa 4021 7661}%
\special{pa 3990 7655}%
\special{pa 3955 7649}%
\special{pa 3918 7645}%
\special{pa 3878 7640}%
\special{pa 3837 7635}%
\special{pa 3797 7628}%
\special{pa 3757 7620}%
\special{pa 3719 7609}%
\special{pa 3683 7595}%
\special{pa 3651 7578}%
\special{pa 3624 7557}%
\special{pa 3600 7533}%
\special{pa 3582 7507}%
\special{pa 3569 7479}%
\special{pa 3561 7450}%
\special{pa 3559 7419}%
\special{pa 3562 7387}%
\special{pa 3567 7355}%
\special{pa 3574 7322}%
\special{pa 3578 7305}%
\special{fp}%
%
\special{pn 8}%
\special{pa 2095 6232}%
\special{pa 4744 6239}%
\special{dt 0.045}%
%
\special{pn 8}%
\special{pa 2743 6786}%
\special{pa 5392 6794}%
\special{dt 0.045}%
%
\special{pn 8}%
\special{pa 3268 7377}%
\special{pa 5752 7377}%
\special{dt 0.045}%
%
\special{pn 8}%
\special{pa 3873 7946}%
\special{pa 6386 7946}%
\special{dt 0.045}%
\put(30.4500,-66.7100){\makebox(0,0)[lb]{$p_0$}}%
\put(14.9000,-71.3200){\makebox(0,0)[lb]{$q_1$}}%
\put(20.3700,-76.9400){\makebox(0,0)[lb]{$q_0$}}%
\put(26.4200,-82.1200){\makebox(0,0)[lb]{$q_{-1}$}}%
\put(23.9000,-61.5300){\makebox(0,0)[lb]{$p_1$}}%
\put(36.0700,-72.1100){\makebox(0,0)[lb]{$p_{-1}$}}%
\put(41.6800,-78.0900){\makebox(0,0)[lb]{$p_{-2}$}}%
%
\special{pn 8}%
\special{pa 2730 6790}%
\special{pa 1930 6780}%
\special{dt 0.045}%
%
\special{pn 8}%
\special{pa 2280 7250}%
\special{pa 2280 6490}%
\special{fp}%
%
\special{pn 8}%
\special{pa 2850 7780}%
\special{pa 2850 6950}%
\special{fp}%
%
\special{pn 8}%
\special{pa 3250 7380}%
\special{pa 2550 7380}%
\special{dt 0.045}%
%
\special{pn 8}%
\special{pa 3850 7940}%
\special{pa 3030 7940}%
\special{dt 0.045}%
\special{pa 3030 7940}%
\special{pa 3030 7940}%
\special{dt 0.045}%
%
\special{pn 8}%
\special{pn 8}%
\special{pa 1910 6780}%
\special{pa 1901 6780}%
\special{fp}%
\special{pa 1862 6781}%
\special{pa 1854 6780}%
\special{fp}%
\special{pa 1814 6777}%
\special{pa 1806 6776}%
\special{fp}%
\special{pa 1768 6765}%
\special{pa 1760 6762}%
\special{fp}%
\special{pa 1724 6746}%
\special{pa 1717 6742}%
\special{fp}%
\special{pa 1682 6722}%
\special{pa 1675 6718}%
\special{fp}%
\special{pa 1640 6701}%
\special{pa 1632 6698}%
\special{fp}%
\special{pa 1593 6690}%
\special{pa 1585 6689}%
\special{fp}%
\special{pa 1546 6690}%
\special{pa 1537 6691}%
\special{fp}%
\special{pa 1498 6698}%
\special{pa 1490 6700}%
\special{fp}%
%
\special{pn 8}%
\special{pn 8}%
\special{pa 2570 7380}%
\special{pa 2562 7379}%
\special{fp}%
\special{pa 2524 7377}%
\special{pa 2516 7376}%
\special{fp}%
\special{pa 2479 7372}%
\special{pa 2471 7371}%
\special{fp}%
\special{pa 2434 7366}%
\special{pa 2426 7365}%
\special{fp}%
\special{pa 2389 7358}%
\special{pa 2381 7357}%
\special{fp}%
\special{pa 2344 7348}%
\special{pa 2336 7346}%
\special{fp}%
\special{pa 2300 7338}%
\special{pa 2292 7336}%
\special{fp}%
\special{pa 2255 7329}%
\special{pa 2247 7328}%
\special{fp}%
\special{pa 2210 7332}%
\special{pa 2203 7335}%
\special{fp}%
\special{pa 2172 7356}%
\special{pa 2166 7361}%
\special{fp}%
\special{pa 2144 7392}%
\special{pa 2141 7399}%
\special{fp}%
\special{pa 2124 7433}%
\special{pa 2121 7441}%
\special{fp}%
\special{pa 2106 7475}%
\special{pa 2103 7483}%
\special{fp}%
\special{pa 2087 7517}%
\special{pa 2083 7523}%
\special{fp}%
\special{pa 2054 7542}%
\special{pa 2047 7540}%
\special{fp}%
\special{pa 2023 7511}%
\special{pa 2019 7504}%
\special{fp}%
\special{pa 2002 7470}%
\special{pa 1999 7463}%
\special{fp}%
\special{pa 1986 7427}%
\special{pa 1984 7420}%
\special{fp}%
\special{pa 1972 7384}%
\special{pa 1970 7376}%
\special{fp}%
\special{pa 1961 7340}%
\special{pa 1959 7332}%
\special{fp}%
\special{pa 1950 7295}%
\special{pa 1949 7288}%
\special{fp}%
\special{pa 1942 7251}%
\special{pa 1941 7243}%
\special{fp}%
\special{pa 1935 7205}%
\special{pa 1934 7197}%
\special{fp}%
\special{pa 1929 7160}%
\special{pa 1928 7152}%
\special{fp}%
\special{pa 1924 7115}%
\special{pa 1923 7107}%
\special{fp}%
\special{pa 1917 7070}%
\special{pa 1915 7062}%
\special{fp}%
\special{pa 1903 7026}%
\special{pa 1900 7019}%
\special{fp}%
\special{pa 1880 6988}%
\special{pa 1875 6982}%
\special{fp}%
\special{pa 1845 6958}%
\special{pa 1839 6954}%
\special{fp}%
\special{pa 1805 6936}%
\special{pa 1798 6933}%
\special{fp}%
\special{pa 1762 6921}%
\special{pa 1755 6919}%
\special{fp}%
\special{pa 1718 6909}%
\special{pa 1710 6907}%
\special{fp}%
\special{pa 1674 6900}%
\special{pa 1666 6899}%
\special{fp}%
\special{pa 1629 6893}%
\special{pa 1620 6892}%
\special{fp}%
\special{pa 1583 6887}%
\special{pa 1575 6886}%
\special{fp}%
\special{pa 1538 6881}%
\special{pa 1530 6880}%
\special{fp}%
%
\special{pn 8}%
\special{pn 8}%
\special{pa 3010 7930}%
\special{pa 3002 7928}%
\special{fp}%
\special{pa 2966 7917}%
\special{pa 2958 7915}%
\special{fp}%
\special{pa 2922 7906}%
\special{pa 2914 7904}%
\special{fp}%
\special{pa 2877 7897}%
\special{pa 2869 7896}%
\special{fp}%
\special{pa 2831 7895}%
\special{pa 2823 7896}%
\special{fp}%
\special{pa 2786 7903}%
\special{pa 2778 7905}%
\special{fp}%
\special{pa 2746 7925}%
\special{pa 2740 7931}%
\special{fp}%
\special{pa 2717 7959}%
\special{pa 2712 7966}%
\special{fp}%
\special{pa 2694 7999}%
\special{pa 2691 8007}%
\special{fp}%
\special{pa 2677 8042}%
\special{pa 2675 8050}%
\special{fp}%
\special{pa 2663 8085}%
\special{pa 2661 8093}%
\special{fp}%
\special{pa 2650 8129}%
\special{pa 2648 8137}%
\special{fp}%
\special{pa 2638 8173}%
\special{pa 2636 8181}%
\special{fp}%
\special{pa 2624 8217}%
\special{pa 2621 8225}%
\special{fp}%
\special{pa 2595 8231}%
\special{pa 2591 8225}%
\special{fp}%
\special{pa 2576 8190}%
\special{pa 2573 8182}%
\special{fp}%
\special{pa 2563 8146}%
\special{pa 2561 8138}%
\special{fp}%
\special{pa 2551 8102}%
\special{pa 2549 8094}%
\special{fp}%
\special{pa 2541 8057}%
\special{pa 2539 8049}%
\special{fp}%
\special{pa 2531 8012}%
\special{pa 2530 8005}%
\special{fp}%
\special{pa 2522 7968}%
\special{pa 2520 7960}%
\special{fp}%
\special{pa 2513 7923}%
\special{pa 2511 7915}%
\special{fp}%
\special{pa 2504 7878}%
\special{pa 2502 7870}%
\special{fp}%
\special{pa 2495 7833}%
\special{pa 2493 7825}%
\special{fp}%
\special{pa 2486 7788}%
\special{pa 2484 7780}%
\special{fp}%
\special{pa 2477 7743}%
\special{pa 2475 7735}%
\special{fp}%
\special{pa 2469 7698}%
\special{pa 2467 7690}%
\special{fp}%
\special{pa 2460 7653}%
\special{pa 2458 7645}%
\special{fp}%
\special{pa 2452 7608}%
\special{pa 2450 7600}%
\special{fp}%
%
\special{pn 8}%
\special{pa 2850 6980}%
\special{pa 2861 6943}%
\special{pa 2871 6906}%
\special{pa 2881 6870}%
\special{pa 2889 6835}%
\special{pa 2896 6800}%
\special{pa 2900 6768}%
\special{pa 2902 6736}%
\special{pa 2901 6707}%
\special{pa 2896 6680}%
\special{pa 2888 6655}%
\special{pa 2875 6633}%
\special{pa 2858 6614}%
\special{pa 2837 6597}%
\special{pa 2813 6582}%
\special{pa 2785 6568}%
\special{pa 2756 6556}%
\special{pa 2724 6545}%
\special{pa 2691 6534}%
\special{pa 2658 6522}%
\special{pa 2624 6511}%
\special{pa 2590 6498}%
\special{pa 2558 6484}%
\special{pa 2526 6469}%
\special{pa 2496 6451}%
\special{pa 2468 6432}%
\special{pa 2443 6410}%
\special{pa 2421 6386}%
\special{pa 2401 6361}%
\special{pa 2385 6334}%
\special{pa 2372 6305}%
\special{pa 2363 6275}%
\special{pa 2356 6244}%
\special{pa 2353 6212}%
\special{pa 2351 6180}%
\special{pa 2350 6147}%
\special{pa 2350 6130}%
\special{fp}%
\end{picture}}%

\caption{}
\end{figure}

{\sc Acknowledgement.} Hearty thanks are due to the referee for careful
reading and many valuable suggestions.

\section{Construction of the diffeomorphism}

\begin{notation}
Denote by $\tau:\R^2\to\R^2$ the translation by $(-1,1)$.
Let
$$\Delta=\{x+y=0\}\subset\R^2.$$
Denote by $\sigma$ the symmetry at $\Delta$: $\sigma(x,y)=(-y,-x)$.
\end{notation}

Notice that $\sigma\tau=\tau\sigma$ and $\sigma^2={\rm id}$. The $C^1$ diffeomorphism $F$ that
we are going to construct will satisfy the following two properties.
\begin{equation}
 \label{3}
\tau F= F\tau, 
\end{equation}
\begin{equation}
 \label{4}
 F^{-1}= \sigma F\sigma.
\end{equation}
Let
$$
P=[-2,\infty)\times[0,1] \ \ \mbox{ and }\ \ 
P'=[0,\infty)\times[0,1].$$
We shall define a surjective diffeomorphism $\phi:P\to P'$ of the form
\begin{equation}
 \label{e1}
\phi(x,y)=(h_y(x),g(y)), 
\end{equation}
where $g:[0,1]\to[0,1]$ is a diffeomorphism with the following properties:
\\[2mm]
(A) $g$ is the time one map of a $C^1$ flow $g^t$ {of the
interval [0,1]},
\\[2mm]
(B) {for any $t\neq0$,} $g^t({y})={y}$ if and only if ${y}\in\{0,\tfrac{1}{4}\}\cup[\tfrac{1}{2},1]$,
\\[2mm]
(C) $g^t\vert_{ \ [0,\tfrac{1}{2}]}$ is symmetric at $\tfrac{1}{4}$, that is
$$
g^t(\tfrac{1}{2}-y)=\tfrac{1}{2}-g^t(y), \ \ \forall y\in[0,\tfrac{1}{2}],$$
(D) $g(y)<y$ for $y\in(0,\tfrac{1}{4})$,
\\[2mm]
(E) $g$ is $C^1$ tangent to the identity at $y=0$,
\\[2mm]
(F) $g$ is affine of slope $e^{\lambda}$ on the
interval $[\tfrac{1}{4}-\delta,\tfrac{1}{4}+\delta]$, where $\lambda$ is some positive
number and  $\delta$ is some {\em small} positive number.
\begin{figure}[h]
{\unitlength 0.1in%
\begin{picture}( 22.3200, 22.4000)( 13.7600,-26.5000)%
%
\special{pn 8}%
\special{pa 1376 2650}%
\special{pa 3608 410}%
\special{fp}%
%
\special{pn 20}%
\special{pa 2488 1546}%
\special{pa 3592 442}%
\special{fp}%
%
\special{pn 20}%
\special{pa 1976 1962}%
\special{pa 1896 2194}%
\special{fp}%
%
\special{pn 20}%
\special{pa 1976 1978}%
\special{pa 1985 1947}%
\special{pa 1995 1916}%
\special{pa 2009 1888}%
\special{pa 2027 1862}%
\special{pa 2048 1838}%
\special{pa 2072 1816}%
\special{pa 2098 1796}%
\special{pa 2125 1777}%
\special{pa 2151 1759}%
\special{pa 2178 1741}%
\special{pa 2205 1724}%
\special{pa 2233 1707}%
\special{pa 2287 1675}%
\special{pa 2315 1659}%
\special{pa 2342 1643}%
\special{pa 2398 1611}%
\special{pa 2425 1594}%
\special{pa 2453 1578}%
\special{pa 2480 1562}%
\special{fp}%
%
\special{pn 20}%
\special{pa 1912 2178}%
\special{pa 1902 2209}%
\special{pa 1890 2239}%
\special{pa 1877 2268}%
\special{pa 1861 2295}%
\special{pa 1841 2319}%
\special{pa 1818 2341}%
\special{pa 1793 2362}%
\special{pa 1766 2381}%
\special{pa 1739 2398}%
\special{pa 1711 2416}%
\special{pa 1655 2450}%
\special{pa 1628 2466}%
\special{pa 1602 2484}%
\special{pa 1575 2501}%
\special{pa 1549 2520}%
\special{pa 1524 2539}%
\special{pa 1500 2559}%
\special{pa 1475 2580}%
\special{pa 1452 2602}%
\special{pa 1428 2623}%
\special{pa 1408 2642}%
\special{fp}%
\put(35.9200,-6.6600){\makebox(0,0)[lb]{$1$}}%
\put(24.8800,-16.7400){\makebox(0,0)[lb]{$1/2$}}%
\put(20.2400,-22.0200){\makebox(0,0)[lb]{$1/4$}}%
\put(14.8000,-27.7800){\makebox(0,0)[lb]{$0$}}%
\end{picture}}%

\caption{The graph of $g$}
\end{figure}
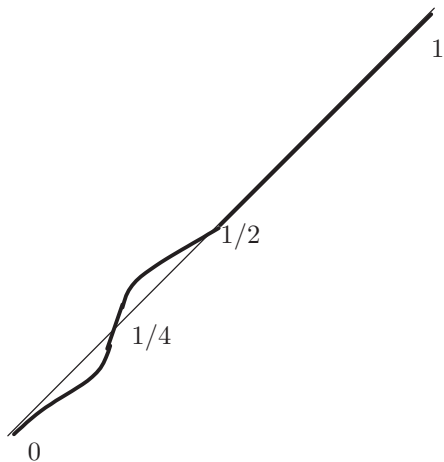

We also assume:
\\[2mm]
(G) $\phi$ sends the rectangle $[-2,-1]\times[0,1]$ onto
 $[0,1]\times
[0,1]$.

\smallskip
The positive number $\lambda$ and a small positive
number $\delta$ will appear in many places. One can show that 
there are such numbers which satisfy all the requirements we pose below.

\begin{remark}
 Notice that the diffeomorphism $\phi$ of form (\ref{e1}) preserves
the horizontal foliations, while it sends the vertical foliation to
itself in the region where $h_y(x)$ does not depend on $y$.
\end{remark}

Let $P_n=\tau^n(P)$ and $P_n'=\tau^n(P')$. (Thus $P=P_0$ and $P'=P_0'$.)
We shall define a diffeomorphism $\phi_n:P_n\to P_n'$ by
$\phi_n=\tau^n\circ \phi\circ\tau^{-n}$ and 
$\Phi:\cup_nP_n\to \cup_n P_n'$ as the union of the $\phi_n$'s.
In order that this defines a homeomorphism, we need the following
condition:
\begin{equation}
 \label{5}
h_1(x-1)=h_0(x)-1, \forall x\geq 1.
\end{equation} 
Of course for $\Phi$ to be a $C^1$ diffeomorphism, we need a bit more.

Define a map $F:\cup_nP_n\to \cup_n P_n'$ by
$$
F=\tau\circ\Phi.$$

Consider a map $F':\sigma(\cup_nP_n)\to\sigma(\cup_nP_n')$ defined by
the conjugation 
$$F'=\sigma\circ F\circ \sigma.
$$
The map $F$ sends the rectangle $[-2,-1]\times[0,1]$ to
$[-1,0]\times[1,2]$,
and reciprocally $F'$ sends $[-1,0]\times[1,2]$ to $[-2,-1]\times[0,1]$.
Routine computation shows that the condition for $F'$ to be
 the inverse of $F$ on these rectangles is the following:
\begin{equation}
 \label{6}
\phi(x,y)=(-g^{-1}(-x-1)+1,g(y)), \ \ \ \forall (x,y)\in[-2,-1]\times[0,1].
\end{equation}
With this condition, we can define a diffeomorphism $F:\R^2\to\R^2$ by
setting it to be equal to $F$ on $\cup_nP_n$ and equal to $(F')^{-1}$
on $\cup_n\sigma(P_n')$. Clearly it satisfies {(\ref{3}) and (\ref{4}).}

\smallskip
Besides (\ref{6}), we assume further conditions on $\phi$:
on $[-1,0]\times[0,\tfrac{1}{2}]$, it is the conjugate of $\phi\vert_{\ [-2,-1]\times[0,\tfrac{1}{2}]}$ by the
translation by (1,0):
\begin{equation}
 \label{7}
\phi(x,y)=(-g^{-1}(-x)+2,g(y)) \ \ \mbox{ on } [-1,0]\times[0,\tfrac{1}{2}].
\end{equation}
{ This condition is helpful
to make the assembled map $\Phi$ to be a $C^1$ diffeomorphism.}
Moreover we assume the following. {
\begin{equation}
\label{8}
\phi(x,y)=(x+2,g(y)) \ \ \ \mbox{on}\ \ \ C,
\end{equation}
where $C$ is the union of the following subsets:
\begin{quote}
 $[-1,\infty)\times[\tfrac{3}{4},1]$, \ \ 
$[0,\infty)\times[\tfrac{1}{2},1]$, \ \  $[0,\infty)\times\{0\}$ and $[0,1]\times[0,1]$.
\end{quote}}
\noindent
 See Figure 3.
\begin{figure}[h]
{\unitlength 0.1in%
\begin{picture}( 42.6400,  8.4000)( 13.7000,-16.3000)%
%
\special{pn 8}%
\special{pa 1394 798}%
\special{pa 5530 798}%
\special{fp}%
%
\special{pn 8}%
\special{pa 1410 1622}%
\special{pa 5634 1614}%
\special{fp}%
%
\special{pn 8}%
\special{pa 1386 798}%
\special{pa 1386 1606}%
\special{fp}%
%
\special{pn 8}%
\special{pa 2258 790}%
\special{pa 2250 1614}%
\special{fp}%
\special{pa 2250 1614}%
\special{pa 2250 1614}%
\special{fp}%
%
\special{pn 8}%
\special{pa 3770 1206}%
\special{pa 3770 1598}%
\special{fp}%
%
\special{pn 8}%
\special{pa 3770 1222}%
\special{pa 5626 1222}%
\special{fp}%
\put(16.9000,-12.6200){\makebox(0,0)[lb]{$A$}}%
\put(24.8200,-14.5400){\makebox(0,0)[lb]{$B$}}%
\put(33.6200,-11.1000){\makebox(0,0)[lb]{$C$}}%
\put(13.7000,-17.5000){\makebox(0,0)[lb]{$-2$}}%
\put(21.7800,-17.4200){\makebox(0,0)[lb]{$-1$}}%
\put(36.9800,-17.2600){\makebox(0,0)[lb]{$1$}}%
\put(29.7000,-17.4000){\makebox(0,0)[lb]{$0$}}%
%
\special{pn 8}%
\special{pa 2250 1000}%
\special{pa 3030 1000}%
\special{fp}%
%
\special{pn 8}%
\special{pa 3020 1000}%
\special{pa 3020 1630}%
\special{fp}%
\put(24.9000,-11.7000){\makebox(0,0)[lb]{$D$}}%
%
\special{pn 8}%
\special{pa 2250 1220}%
\special{pa 3020 1220}%
\special{fp}%
\put(41.6000,-14.7000){\makebox(0,0)[lb]{$G$}}%
\end{picture}}%

\caption{Rectangles $A$ and $B$ are mapped by $\phi$ 
{onto the rectangles two unit right to them, by a product map (\ref{6})
 and
(\ref{7})}. On the other hand,
 $\phi(x,y)=(x+2,g(y))$ on $C$. }
\end{figure}
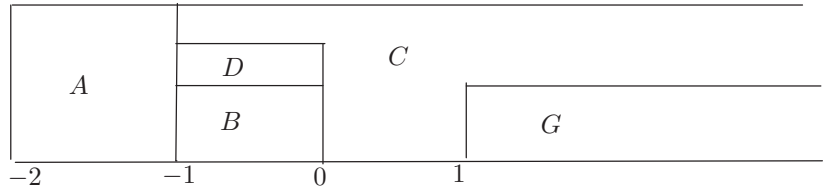

{ The map $\phi$ is already determined on
the boundary of
 $D=[-1,0]\times[\tfrac{1}{2},\tfrac{1}{4}]$.
 On $D$, $\phi$ is to be
any extension of it of the form (\ref{e1}).}
  Notice that
 the map $\phi$ defined by (\ref{6}), (\ref{7}) and (\ref{8})
satisfies the condition (\ref{5}).

The foliation $\FF^u$ is defined to be the image by the iterates of $F$ of
 the vertical foliation on
$\cup_n\sigma(P_n)$.
Conversely $\FF^s$ is to be the image by the iterates of $F^{-1}$ of the
horizontal foliation on $\cup_n P_n$.
 More concretely on $P$, $\FF^s$ is the horizontal foliation, while
$\FF^u$ is the image by the iterates of
$\phi$ of the vertical foliation on $[-2,0]\times[0,1]$.
Since the product map sends the vertical foliation to the vertical
 foliation, we have:
\begin{lemma}
 \label{lemma}
The foliation $\FF^u$ is vertical on $[-2,\infty)\times[\tfrac{3}{4},1]$
and also on $[-2,3]\times[0,\tfrac{1}{2}]$.
\end{lemma}

\section{More conditions on the map $\phi$}

{In this section,
we shall define a map $\phi$ on the region $G$ in Figure 3.
For this, we 
first define a Reeb component $R$ of the foliation $\FF^u$ in $P$ as in
Figure 5.}
Let us define its boundary $\partial R$ to be
the graph of a function $\theta:[0,\tfrac{1}{2}]\to[4,\infty)$
symmetric at $y=\tfrac{1}{4}$. By the symmetry, we need to define $\theta$ only on $[0,\tfrac{1}{4}]$.
Recall that the map $g:[0,\tfrac{1}{2}]\to[0,\tfrac{1}{2}]$ is the time one map of the
flow
$g^t$. If we put $y(t)=g^t(1/8)$, it is
monotone decreasing and satisfies {
$\lim_{t\to-\infty}y(t)=\tfrac{1}{4}$ and
$\lim_{t\to\infty}y(t)=0$.}
 Let us define first of all
a curve
$x(t)\in[4,\infty)$, $t\in\R$,  and then
a function $\theta$ by $\theta(y(t))=x(t)$.
The conditions for $x(t)$ are the following:
\\[2mm]
(H) $x(t)=4$ for $t<t_0$ for some $t_0<0$, equivalently, $\theta(y)=4$
if $y$ is $\delta$-near to $\tfrac{1}{4}$, where $\delta=\tfrac{1}{4}-y(t_0)>0$ is some small number.
\\[2mm]
(I) $x'(t)\in[0,2)$ and $x'(t)$ is {strictly} monotone increasing for $t>t_0$ and
{$\lim_{t\to\infty}x'(t)=2$.}
\\[2mm]
Thus { $x(t)$ itself is monotone increasing. Moreover,} we have $x(t+1)<x(t)+2$ and its difference tends to 0 monotonically.
Define the Reeb component $R$ 
by
$$
R=\{x\geq\theta(y), 0<y<\tfrac{1}{2}\}.$$
We have
$$
\partial R={\rm Graph}(\theta).
$$
{Conditions (H) and (I) imply that $\partial R$ is vertical on
the region $\abs{y-\tfrac{1}{4}}<\delta$ and is strictly convex leftward
 outside this region.}

Next we shall define the diffeomorphism 
$$\phi:([0,\infty)\times[0,\tfrac{1}{2}])\setminus R \to
([2,\infty)\times[0,\tfrac{1}{2}])\setminus R. $$
Again $\phi$ is to be symmetric with respect
 to the line $y=\tfrac{1}{4}$, and we shall define it only 
on $([0,\infty)\times[0,\tfrac{1}{4}])\setminus R$.

\smallskip
$\bullet$ On $\partial R\cap([0,\infty)\times[0,\tfrac{1}{4}])$, $\phi$ is defined by
$\phi(x(t),y(t))=(x(t+1),y(t+1))$. 

\smallskip
$\bullet$ On 
$([0,\infty)\times[0,\tfrac{1}{4}])\setminus {\Int}(R)$, $\phi$ maps the interval 
$[0,x(t)]\times\{y(t)\}$ { to the interval} $[2,x(t+1)]\times \{y(t+1)\}$ by the formula
$$\phi(x,y(t))=(h_{y(t)}(x),y(t+1)),$$ where the diffeomorphism 
$$h_{y(t)}:[0,x(t)]\to[2,x(t+1)]$$
satisfies:
\begin{equation}
 \label{9}
h_{y(t)}(x)=x+2 \ \ \ \mbox{ if }\ \ \ 0\leq x\leq x(t)/4, 
\end{equation}
\begin{equation}
 \label{10}
\Lambda=\{(x,y(t))\mid h_{y(t)}'(x)=e^{-\lambda}\} \ \ 
\mbox{ is a  neighbourhood of $\partial R$ in $P\setminus{\rm Int}(R)$.}
\end{equation}
{Recall that $\lambda>0$ is a constant which appeared in
 condition  (F) on $g$.}
\begin{equation}
\label{11}
h_{y(t)}'(x)\leq 1, \ \ \forall x\in[0,x(t)].
\end{equation}
\begin{equation}
 \label{12}
h_{y(t)}\ \ \mbox{does not depend on $y(t)$ if }\ \tfrac{1}{4}-\delta\leq g(y(t))=y(t+1)\leq \tfrac{1}{4}.
\end{equation}

The following lemma is a restatement of (\ref{9}).
See Figure 4.

\begin{lemma}
 \label{lemma4}
On the region $\{0<x<\theta(y)/4,\ 0<y<\tfrac{1}{2}\}$, we have
 $$\phi(x,y)=(x+2,g(y)).$$ \qed
\end{lemma}
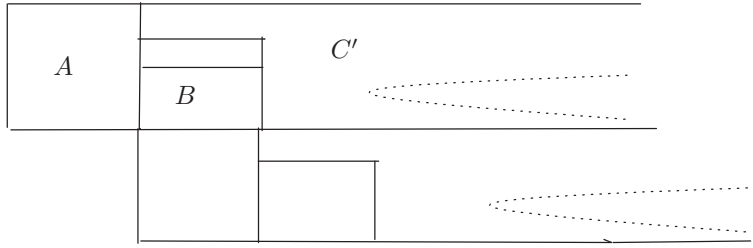
\begin{figure}[h]
{\unitlength 0.1in%
\begin{picture}( 38.9100, 12.6000)( 14.0600,-21.3000)%
%
\special{pn 8}%
\special{pa 1412 876}%
\special{pa 4721 876}%
\special{fp}%
%
\special{pn 8}%
\special{pa 1425 1536}%
\special{pa 4804 1529}%
\special{fp}%
%
\special{pn 8}%
\special{pa 1406 876}%
\special{pa 1406 1523}%
\special{fp}%
%
\special{pn 8}%
\special{pa 2104 870}%
\special{pa 2097 1529}%
\special{fp}%
\special{pa 2097 1529}%
\special{pa 2097 1529}%
\special{fp}%
%
\special{pn 8}%
\special{pa 2116 1209}%
\special{pa 2744 1209}%
\special{fp}%
\put(16.4900,-12.4800){\makebox(0,0)[lb]{$A$}}%
\put(22.8300,-14.0100){\makebox(0,0)[lb]{$B$}}%
%
\special{pn 8}%
\special{pa 2105 2118}%
\special{pa 4569 2126}%
\special{fp}%
\special{pa 4569 2126}%
\special{pa 4529 2110}%
\special{fp}%
%
\special{pn 8}%
\special{pa 2721 1526}%
\special{pa 2721 2126}%
\special{fp}%
%
\special{pn 8}%
\special{pa 4577 2126}%
\special{pa 5297 2118}%
\special{fp}%
%
\special{pn 8}%
\special{pa 2740 1050}%
\special{pa 2740 1540}%
\special{fp}%
%
\special{pn 8}%
\special{pa 2090 1060}%
\special{pa 2750 1060}%
\special{fp}%
\special{pa 2750 1060}%
\special{pa 2750 1060}%
\special{fp}%
%
\special{pn 8}%
\special{pa 2720 1700}%
\special{pa 3350 1700}%
\special{fp}%
%
\special{pn 8}%
\special{pa 3330 1700}%
\special{pa 3330 2110}%
\special{fp}%
%
\special{pn 8}%
\special{pn 8}%
\special{pa 4650 1250}%
\special{pa 4642 1250}%
\special{fp}%
\special{pa 4605 1252}%
\special{pa 4597 1252}%
\special{fp}%
\special{pa 4560 1254}%
\special{pa 4552 1254}%
\special{fp}%
\special{pa 4515 1255}%
\special{pa 4507 1255}%
\special{fp}%
\special{pa 4470 1257}%
\special{pa 4462 1257}%
\special{fp}%
\special{pa 4425 1258}%
\special{pa 4417 1259}%
\special{fp}%
\special{pa 4380 1260}%
\special{pa 4372 1261}%
\special{fp}%
\special{pa 4335 1262}%
\special{pa 4327 1263}%
\special{fp}%
\special{pa 4290 1264}%
\special{pa 4282 1265}%
\special{fp}%
\special{pa 4245 1266}%
\special{pa 4237 1266}%
\special{fp}%
\special{pa 4200 1268}%
\special{pa 4192 1268}%
\special{fp}%
\special{pa 4155 1269}%
\special{pa 4147 1269}%
\special{fp}%
\special{pa 4110 1271}%
\special{pa 4102 1271}%
\special{fp}%
\special{pa 4065 1273}%
\special{pa 4057 1274}%
\special{fp}%
\special{pa 4020 1275}%
\special{pa 4012 1276}%
\special{fp}%
\special{pa 3975 1277}%
\special{pa 3968 1278}%
\special{fp}%
\special{pa 3931 1279}%
\special{pa 3923 1280}%
\special{fp}%
\special{pa 3886 1282}%
\special{pa 3878 1282}%
\special{fp}%
\special{pa 3841 1284}%
\special{pa 3833 1284}%
\special{fp}%
\special{pa 3796 1286}%
\special{pa 3788 1287}%
\special{fp}%
\special{pa 3751 1289}%
\special{pa 3743 1289}%
\special{fp}%
\special{pa 3706 1291}%
\special{pa 3698 1292}%
\special{fp}%
\special{pa 3661 1294}%
\special{pa 3653 1294}%
\special{fp}%
\special{pa 3616 1297}%
\special{pa 3608 1297}%
\special{fp}%
\special{pa 3571 1299}%
\special{pa 3563 1300}%
\special{fp}%
\special{pa 3526 1303}%
\special{pa 3518 1304}%
\special{fp}%
\special{pa 3481 1307}%
\special{pa 3473 1307}%
\special{fp}%
\special{pa 3437 1311}%
\special{pa 3429 1312}%
\special{fp}%
\special{pa 3392 1316}%
\special{pa 3384 1317}%
\special{fp}%
\special{pa 3347 1323}%
\special{pa 3340 1324}%
\special{fp}%
\special{pa 3304 1335}%
\special{pa 3301 1341}%
\special{fp}%
\special{pa 3333 1356}%
\special{pa 3341 1357}%
\special{fp}%
\special{pa 3377 1365}%
\special{pa 3385 1366}%
\special{fp}%
\special{pa 3422 1372}%
\special{pa 3430 1373}%
\special{fp}%
\special{pa 3466 1378}%
\special{pa 3474 1379}%
\special{fp}%
\special{pa 3511 1383}%
\special{pa 3519 1384}%
\special{fp}%
\special{pa 3556 1388}%
\special{pa 3564 1390}%
\special{fp}%
\special{pa 3601 1393}%
\special{pa 3609 1394}%
\special{fp}%
\special{pa 3645 1398}%
\special{pa 3653 1399}%
\special{fp}%
\special{pa 3690 1403}%
\special{pa 3698 1403}%
\special{fp}%
\special{pa 3735 1407}%
\special{pa 3743 1408}%
\special{fp}%
\special{pa 3780 1411}%
\special{pa 3788 1412}%
\special{fp}%
\special{pa 3825 1415}%
\special{pa 3833 1415}%
\special{fp}%
\special{pa 3869 1419}%
\special{pa 3877 1420}%
\special{fp}%
\special{pa 3914 1423}%
\special{pa 3922 1423}%
\special{fp}%
\special{pa 3959 1427}%
\special{pa 3967 1427}%
\special{fp}%
\special{pa 4004 1430}%
\special{pa 4012 1431}%
\special{fp}%
\special{pa 4049 1434}%
\special{pa 4057 1434}%
\special{fp}%
\special{pa 4094 1438}%
\special{pa 4102 1438}%
\special{fp}%
\special{pa 4138 1442}%
\special{pa 4146 1442}%
\special{fp}%
\special{pa 4183 1445}%
\special{pa 4191 1445}%
\special{fp}%
\special{pa 4228 1449}%
\special{pa 4236 1449}%
\special{fp}%
\special{pa 4273 1452}%
\special{pa 4281 1452}%
\special{fp}%
\special{pa 4318 1455}%
\special{pa 4326 1456}%
\special{fp}%
\special{pa 4363 1459}%
\special{pa 4371 1460}%
\special{fp}%
\special{pa 4408 1463}%
\special{pa 4416 1463}%
\special{fp}%
\special{pa 4453 1466}%
\special{pa 4461 1466}%
\special{fp}%
\special{pa 4497 1469}%
\special{pa 4505 1470}%
\special{fp}%
\special{pa 4542 1472}%
\special{pa 4550 1473}%
\special{fp}%
\special{pa 4587 1476}%
\special{pa 4595 1477}%
\special{fp}%
\special{pa 4632 1480}%
\special{pa 4640 1480}%
\special{fp}%
%
\special{pn 8}%
\special{pn 8}%
\special{pa 5280 1840}%
\special{pa 5272 1840}%
\special{fp}%
\special{pa 5235 1842}%
\special{pa 5227 1842}%
\special{fp}%
\special{pa 5190 1844}%
\special{pa 5182 1844}%
\special{fp}%
\special{pa 5145 1845}%
\special{pa 5137 1845}%
\special{fp}%
\special{pa 5100 1847}%
\special{pa 5092 1847}%
\special{fp}%
\special{pa 5055 1848}%
\special{pa 5047 1849}%
\special{fp}%
\special{pa 5010 1850}%
\special{pa 5002 1851}%
\special{fp}%
\special{pa 4965 1852}%
\special{pa 4957 1853}%
\special{fp}%
\special{pa 4920 1854}%
\special{pa 4912 1855}%
\special{fp}%
\special{pa 4875 1856}%
\special{pa 4867 1856}%
\special{fp}%
\special{pa 4830 1858}%
\special{pa 4822 1858}%
\special{fp}%
\special{pa 4785 1859}%
\special{pa 4777 1859}%
\special{fp}%
\special{pa 4740 1861}%
\special{pa 4732 1861}%
\special{fp}%
\special{pa 4695 1863}%
\special{pa 4687 1864}%
\special{fp}%
\special{pa 4650 1865}%
\special{pa 4642 1866}%
\special{fp}%
\special{pa 4605 1867}%
\special{pa 4598 1868}%
\special{fp}%
\special{pa 4561 1869}%
\special{pa 4553 1870}%
\special{fp}%
\special{pa 4516 1872}%
\special{pa 4508 1872}%
\special{fp}%
\special{pa 4471 1874}%
\special{pa 4463 1874}%
\special{fp}%
\special{pa 4426 1876}%
\special{pa 4418 1877}%
\special{fp}%
\special{pa 4381 1879}%
\special{pa 4373 1879}%
\special{fp}%
\special{pa 4336 1881}%
\special{pa 4328 1882}%
\special{fp}%
\special{pa 4291 1884}%
\special{pa 4283 1884}%
\special{fp}%
\special{pa 4246 1887}%
\special{pa 4238 1887}%
\special{fp}%
\special{pa 4201 1889}%
\special{pa 4193 1890}%
\special{fp}%
\special{pa 4156 1893}%
\special{pa 4148 1894}%
\special{fp}%
\special{pa 4111 1897}%
\special{pa 4103 1897}%
\special{fp}%
\special{pa 4067 1901}%
\special{pa 4059 1902}%
\special{fp}%
\special{pa 4022 1906}%
\special{pa 4014 1907}%
\special{fp}%
\special{pa 3977 1913}%
\special{pa 3970 1914}%
\special{fp}%
\special{pa 3934 1925}%
\special{pa 3931 1931}%
\special{fp}%
\special{pa 3963 1946}%
\special{pa 3971 1947}%
\special{fp}%
\special{pa 4007 1955}%
\special{pa 4015 1956}%
\special{fp}%
\special{pa 4052 1962}%
\special{pa 4060 1963}%
\special{fp}%
\special{pa 4096 1968}%
\special{pa 4104 1969}%
\special{fp}%
\special{pa 4141 1973}%
\special{pa 4149 1974}%
\special{fp}%
\special{pa 4186 1978}%
\special{pa 4194 1980}%
\special{fp}%
\special{pa 4231 1983}%
\special{pa 4239 1984}%
\special{fp}%
\special{pa 4275 1988}%
\special{pa 4283 1989}%
\special{fp}%
\special{pa 4320 1993}%
\special{pa 4328 1993}%
\special{fp}%
\special{pa 4365 1997}%
\special{pa 4373 1998}%
\special{fp}%
\special{pa 4410 2001}%
\special{pa 4418 2002}%
\special{fp}%
\special{pa 4455 2005}%
\special{pa 4463 2005}%
\special{fp}%
\special{pa 4499 2009}%
\special{pa 4507 2010}%
\special{fp}%
\special{pa 4544 2013}%
\special{pa 4552 2013}%
\special{fp}%
\special{pa 4589 2017}%
\special{pa 4597 2017}%
\special{fp}%
\special{pa 4634 2020}%
\special{pa 4642 2021}%
\special{fp}%
\special{pa 4679 2024}%
\special{pa 4687 2024}%
\special{fp}%
\special{pa 4724 2028}%
\special{pa 4732 2028}%
\special{fp}%
\special{pa 4768 2032}%
\special{pa 4776 2032}%
\special{fp}%
\special{pa 4813 2035}%
\special{pa 4821 2035}%
\special{fp}%
\special{pa 4858 2039}%
\special{pa 4866 2039}%
\special{fp}%
\special{pa 4903 2042}%
\special{pa 4911 2042}%
\special{fp}%
\special{pa 4948 2045}%
\special{pa 4956 2046}%
\special{fp}%
\special{pa 4993 2049}%
\special{pa 5001 2050}%
\special{fp}%
\special{pa 5038 2053}%
\special{pa 5046 2053}%
\special{fp}%
\special{pa 5083 2056}%
\special{pa 5091 2056}%
\special{fp}%
\special{pa 5127 2059}%
\special{pa 5135 2060}%
\special{fp}%
\special{pa 5172 2062}%
\special{pa 5180 2063}%
\special{fp}%
\special{pa 5217 2066}%
\special{pa 5225 2067}%
\special{fp}%
\special{pa 5262 2070}%
\special{pa 5270 2070}%
\special{fp}%
%
\special{pn 8}%
\special{pa 2090 1530}%
\special{pa 2090 2130}%
\special{fp}%
\put(31.0000,-11.6000){\makebox(0,0)[lb]{$C'$}}%
\end{picture}}%

\caption{
The dotted curve indicates $\{x=\theta(y)/4\}$.
$\phi(x,y)=(x+2,g(y))$ { on the region $C'$ which is outside of the
dotted curve.}}
\end{figure}

This lemma, together with the fact that
$g'(0)=g'(1)=1$, shows that the assembled map
$\Phi:\cup_nP_n\to\cup_nP'_n$ is actually
a $C^1$ diffeomorphism. Denote the Euclidean norm on $\R^2$ by $\abs{\cdot}$.

\begin{corollary}
 \label{corollary}
The tangent bundle $T\FF^u$ of the foliation $\FF^u$ is vertical in
a neighbouhood of $[0,\infty)\times\{0,1\}$ and if $v\in T\FF^u_p$, 
$p\in [0,\infty)\times\{0,1\}$,
then $\abs{D\phi(v)}=\abs{v}$.
\end{corollary}

\bd The first assertion follows from Lemmas \ref{lemma} and
 \ref{lemma4},
 while the last from $g'(0)=g'(1)=1$.  \qed

\medskip

So far we have defined
the diffeomorphism $\phi$, whence the foliation $\FF^u$,
 except in the interior of the Reeb
component $R$. On $R=\{ x\geq\theta(y),\ 0<y<\tfrac{1}{2}\}$, define the
foliation $\FF^u$ by the horizontal translation of the boundary curve
$\partial R$. See Figure 5. Let $L=[4,\infty)\times\{\tfrac{1}{4}\}$ be
the center { ray} of $R$.
The two transverse foliations $\FF^u$ and $\FF^s$ define a product
structure on $R$:
$$R\approx\partial R\times L.$$ 
 We have already defined the map $\phi$ on $\partial R$.
Let us define it on $L$ to be the contraction of ratio $e^{-\lambda}$
centered at $(\tfrac{1}{4},4)$.
Finally define the map $\phi:R\to R$ as the product of these two maps.
By virtue of (\ref{10}), $\phi:P\to P'$ is a $C^1$
diffeomorphism. Recall that it has the form
$\phi(x,y)=(h_y(x),g(y))$. 

\begin{lemma}
 \label{lemma2}
(1) $h_y'(x)\in (0,1]$ for any $(x,y)\in P'.$

(2) There is a neighbourhood $N$ of the point $(4,\tfrac{1}{4})$  such that if
 $(x,y)\in N\cup R$, then $h_y'(x)=e^{-\lambda}$.

(3) {Moreover one can choose $N$ of (2) large enough so that 
 if $(x,y)\in P\setminus(N\cup R)$, then
$h_y(x)-x>\alpha$ for some fixed $\alpha>0$.}
\end{lemma}

\bd {
(1) follows from (\ref{11}) and the construction on $R$. For (2), one can choose $N$ to be any neighbourhood of
 $(4,\tfrac{1}{4})$ in $\Lambda\cup R$, where $\Lambda$ is a set given
 by (\ref{10}).} Let us show (3).
 Conditions (H) and (I) imply that
$x(t)$ is strictly monotone increasing if $x(t)>4$.
This, together with (\ref{9}) and (\ref{11}), shows that the set 
$$K=\{(x,y)\in P\setminus \Int(R)\mid h_y(x)\leq x\}$$
coincides with a compact interval
$$I=\{(x(t),y(t))\mid {x(t+1)=x(t)=4}\}$$
of $\partial R$.  One can choose a neighbourhood $N$ of $I$
contained in the set $\Lambda$,
and set
$$
{\alpha}=\min\{h_y(x)-x\mid(x,y)\in P\setminus(N\cup R)\}.$$ \qed

\begin{figure}[h]
{\unitlength 0.1in%
\begin{picture}( 43.5400, 17.4300)(  4.8000,-22.7700)%
%
\special{pn 13}%
\special{pa 695 595}%
\special{pa 4811 595}%
\special{dt 0.045}%
%
\special{pn 13}%
\special{pa 687 2230}%
\special{pa 4804 2223}%
\special{dt 0.045}%
%
\special{pn 8}%
\special{pa 2699 1302}%
\special{pa 2699 1501}%
\special{fp}%
\special{pa 2699 1501}%
\special{pa 2699 1501}%
\special{fp}%
%
\special{pn 8}%
\special{pa 3130 1317}%
\special{pa 3130 1501}%
\special{fp}%
%
\special{pn 8}%
\special{pa 3514 1325}%
\special{pa 3514 1494}%
\special{fp}%
%
\special{pn 8}%
\special{pa 2699 1325}%
\special{pa 2704 1293}%
\special{pa 2710 1262}%
\special{pa 2717 1231}%
\special{pa 2726 1202}%
\special{pa 2736 1173}%
\special{pa 2747 1146}%
\special{pa 2760 1120}%
\special{pa 2774 1095}%
\special{pa 2789 1071}%
\special{pa 2806 1048}%
\special{pa 2824 1026}%
\special{pa 2843 1004}%
\special{pa 2863 984}%
\special{pa 2907 946}%
\special{pa 2931 929}%
\special{pa 2956 912}%
\special{pa 2982 896}%
\special{pa 3008 881}%
\special{pa 3036 866}%
\special{pa 3065 852}%
\special{pa 3094 839}%
\special{pa 3124 827}%
\special{pa 3155 815}%
\special{pa 3187 804}%
\special{pa 3220 793}%
\special{pa 3253 783}%
\special{pa 3287 773}%
\special{pa 3321 764}%
\special{pa 3357 755}%
\special{pa 3392 747}%
\special{pa 3429 739}%
\special{pa 3465 731}%
\special{pa 3502 724}%
\special{pa 3540 717}%
\special{pa 3616 705}%
\special{pa 3733 687}%
\special{pa 3772 682}%
\special{pa 3852 672}%
\special{pa 3891 666}%
\special{pa 3971 656}%
\special{pa 3974 656}%
\special{fp}%
%
\special{pn 8}%
\special{pa 2699 1486}%
\special{pa 2710 1516}%
\special{pa 2721 1547}%
\special{pa 2732 1576}%
\special{pa 2745 1604}%
\special{pa 2775 1658}%
\special{pa 2791 1683}%
\special{pa 2808 1707}%
\special{pa 2846 1753}%
\special{pa 2866 1774}%
\special{pa 2888 1795}%
\special{pa 2910 1814}%
\special{pa 2933 1833}%
\special{pa 2957 1852}%
\special{pa 2981 1869}%
\special{pa 3033 1901}%
\special{pa 3060 1917}%
\special{pa 3088 1931}%
\special{pa 3117 1945}%
\special{pa 3146 1958}%
\special{pa 3176 1971}%
\special{pa 3206 1983}%
\special{pa 3237 1994}%
\special{pa 3301 2016}%
\special{pa 3334 2026}%
\special{pa 3367 2035}%
\special{pa 3435 2053}%
\special{pa 3470 2061}%
\special{pa 3504 2069}%
\special{pa 3540 2077}%
\special{pa 3575 2084}%
\special{pa 3611 2091}%
\special{pa 3648 2098}%
\special{pa 3684 2104}%
\special{pa 3721 2110}%
\special{pa 3758 2117}%
\special{pa 3795 2122}%
\special{pa 3869 2134}%
\special{pa 3907 2140}%
\special{pa 3944 2145}%
\special{pa 3982 2151}%
\special{pa 4005 2154}%
\special{fp}%
%
\special{pn 8}%
\special{pa 3130 1340}%
\special{pa 3145 1244}%
\special{pa 3152 1214}%
\special{pa 3160 1184}%
\special{pa 3170 1156}%
\special{pa 3181 1129}%
\special{pa 3194 1103}%
\special{pa 3207 1078}%
\special{pa 3223 1054}%
\special{pa 3239 1032}%
\special{pa 3257 1010}%
\special{pa 3276 989}%
\special{pa 3296 969}%
\special{pa 3317 951}%
\special{pa 3340 933}%
\special{pa 3363 916}%
\special{pa 3388 899}%
\special{pa 3440 869}%
\special{pa 3468 856}%
\special{pa 3496 842}%
\special{pa 3556 818}%
\special{pa 3587 807}%
\special{pa 3619 797}%
\special{pa 3652 787}%
\special{pa 3685 778}%
\special{pa 3719 769}%
\special{pa 3789 753}%
\special{pa 3861 739}%
\special{pa 3898 733}%
\special{pa 3936 727}%
\special{pa 3973 721}%
\special{pa 4012 716}%
\special{pa 4050 711}%
\special{pa 4089 707}%
\special{pa 4129 702}%
\special{pa 4168 698}%
\special{pa 4288 686}%
\special{pa 4328 683}%
\special{pa 4369 679}%
\special{pa 4409 675}%
\special{pa 4450 672}%
\special{fp}%
%
\special{pn 8}%
\special{pa 3130 1470}%
\special{pa 3154 1530}%
\special{pa 3167 1559}%
\special{pa 3195 1615}%
\special{pa 3211 1641}%
\special{pa 3245 1691}%
\special{pa 3283 1737}%
\special{pa 3325 1779}%
\special{pa 3347 1799}%
\special{pa 3370 1818}%
\special{pa 3394 1836}%
\special{pa 3419 1854}%
\special{pa 3444 1870}%
\special{pa 3470 1887}%
\special{pa 3524 1917}%
\special{pa 3552 1930}%
\special{pa 3581 1944}%
\special{pa 3610 1957}%
\special{pa 3640 1969}%
\special{pa 3733 2002}%
\special{pa 3766 2012}%
\special{pa 3798 2022}%
\special{pa 3831 2031}%
\special{pa 3865 2039}%
\special{pa 3898 2048}%
\special{pa 3933 2056}%
\special{pa 3967 2063}%
\special{pa 4037 2077}%
\special{pa 4073 2084}%
\special{pa 4109 2090}%
\special{pa 4145 2097}%
\special{pa 4181 2103}%
\special{pa 4217 2108}%
\special{pa 4254 2114}%
\special{pa 4291 2119}%
\special{pa 4327 2125}%
\special{pa 4438 2140}%
\special{pa 4476 2145}%
\special{pa 4482 2146}%
\special{fp}%
%
\special{pn 8}%
\special{pa 3514 1347}%
\special{pa 3522 1283}%
\special{pa 3527 1252}%
\special{pa 3533 1221}%
\special{pa 3541 1192}%
\special{pa 3551 1164}%
\special{pa 3562 1138}%
\special{pa 3574 1112}%
\special{pa 3588 1087}%
\special{pa 3603 1064}%
\special{pa 3620 1041}%
\special{pa 3638 1020}%
\special{pa 3657 999}%
\special{pa 3678 979}%
\special{pa 3699 961}%
\special{pa 3722 943}%
\special{pa 3746 926}%
\special{pa 3771 910}%
\special{pa 3797 894}%
\special{pa 3825 879}%
\special{pa 3853 865}%
\special{pa 3882 852}%
\special{pa 3942 828}%
\special{pa 3974 816}%
\special{pa 4040 796}%
\special{pa 4074 786}%
\special{pa 4108 777}%
\special{pa 4144 768}%
\special{pa 4180 760}%
\special{pa 4216 753}%
\special{pa 4253 745}%
\special{pa 4290 738}%
\special{pa 4328 732}%
\special{pa 4484 708}%
\special{pa 4644 688}%
\special{pa 4685 684}%
\special{pa 4725 679}%
\special{pa 4766 675}%
\special{pa 4789 672}%
\special{fp}%
%
\special{pn 8}%
\special{pa 3514 1470}%
\special{pa 3538 1530}%
\special{pa 3551 1559}%
\special{pa 3579 1615}%
\special{pa 3595 1641}%
\special{pa 3612 1666}%
\special{pa 3630 1690}%
\special{pa 3668 1736}%
\special{pa 3689 1757}%
\special{pa 3711 1778}%
\special{pa 3733 1797}%
\special{pa 3756 1816}%
\special{pa 3780 1834}%
\special{pa 3805 1852}%
\special{pa 3857 1884}%
\special{pa 3884 1899}%
\special{pa 3912 1913}%
\special{pa 3941 1927}%
\special{pa 3970 1940}%
\special{pa 4000 1953}%
\special{pa 4030 1965}%
\special{pa 4092 1987}%
\special{pa 4124 1997}%
\special{pa 4157 2007}%
\special{pa 4223 2025}%
\special{pa 4257 2034}%
\special{pa 4291 2042}%
\special{pa 4325 2049}%
\special{pa 4360 2057}%
\special{pa 4395 2064}%
\special{pa 4431 2071}%
\special{pa 4467 2077}%
\special{pa 4503 2084}%
\special{pa 4575 2096}%
\special{pa 4612 2101}%
\special{pa 4648 2107}%
\special{pa 4685 2113}%
\special{pa 4759 2123}%
\special{pa 4796 2129}%
\special{pa 4811 2131}%
\special{fp}%
%
\special{pn 20}%
\special{pa 2231 1332}%
\special{pa 749 1332}%
\special{fp}%
%
\special{pn 20}%
\special{pa 2246 1532}%
\special{pa 772 1532}%
\special{fp}%
%
\special{pn 20}%
\special{pa 2008 1332}%
\special{pa 2008 1524}%
\special{fp}%
\put(22.5400,-24.0700){\makebox(0,0)[lb]{$4$}}%
\put(4.8000,-6.6400){\makebox(0,0)[lb]{$1/2$}}%
\put(5.5700,-14.7000){\makebox(0,0)[lb]{$1/4$}}%
\put(5.1800,-22.3000){\makebox(0,0)[lb]{$0$}}%
%
\special{pn 13}%
\special{pa 1962 826}%
\special{pa 4781 833}%
\special{dt 0.045}%
%
\special{pn 13}%
\special{pa 1978 1778}%
\special{pa 4804 1770}%
\special{dt 0.045}%
\special{pa 4804 1770}%
\special{pa 4804 1770}%
\special{dt 0.045}%
%
\special{pn 13}%
\special{pa 1993 2023}%
\special{pa 4834 2023}%
\special{dt 0.045}%
\special{pa 4834 2023}%
\special{pa 4834 2023}%
\special{dt 0.045}%
\put(11.6200,-16.4900){\makebox(0,0)[lb]{$E$}}%
\put(40.9000,-13.6100){\makebox(0,0)[lb]{$L$}}%
\put(24.8600,-7.5600){\makebox(0,0)[lb]{$\partial R$}}%
%
\special{pn 20}%
\special{pa 2224 1328}%
\special{pa 2224 1536}%
\special{fp}%
%
\special{pn 20}%
\special{pa 2216 1344}%
\special{pa 2216 1311}%
\special{pa 2218 1279}%
\special{pa 2224 1248}%
\special{pa 2236 1220}%
\special{pa 2253 1193}%
\special{pa 2274 1167}%
\special{pa 2297 1143}%
\special{pa 2320 1120}%
\special{pa 2370 1076}%
\special{pa 2395 1055}%
\special{pa 2420 1035}%
\special{pa 2472 997}%
\special{pa 2499 980}%
\special{pa 2526 962}%
\special{pa 2553 946}%
\special{pa 2609 914}%
\special{pa 2665 886}%
\special{pa 2694 872}%
\special{pa 2752 846}%
\special{pa 2782 834}%
\special{pa 2812 823}%
\special{pa 2842 811}%
\special{pa 2872 801}%
\special{pa 2965 771}%
\special{pa 3027 753}%
\special{pa 3058 745}%
\special{pa 3154 721}%
\special{pa 3282 693}%
\special{pa 3315 686}%
\special{pa 3347 680}%
\special{pa 3380 673}%
\special{pa 3412 667}%
\special{pa 3511 649}%
\special{pa 3543 643}%
\special{pa 3560 640}%
\special{fp}%
%
\special{pn 20}%
\special{pa 2224 1536}%
\special{pa 2234 1567}%
\special{pa 2246 1596}%
\special{pa 2261 1624}%
\special{pa 2280 1650}%
\special{pa 2300 1675}%
\special{pa 2322 1699}%
\special{pa 2344 1722}%
\special{pa 2367 1744}%
\special{pa 2391 1765}%
\special{pa 2439 1805}%
\special{pa 2464 1824}%
\special{pa 2490 1842}%
\special{pa 2516 1859}%
\special{pa 2543 1876}%
\special{pa 2570 1892}%
\special{pa 2597 1907}%
\special{pa 2625 1922}%
\special{pa 2653 1936}%
\special{pa 2740 1975}%
\special{pa 2770 1986}%
\special{pa 2800 1998}%
\special{pa 2831 2009}%
\special{pa 2861 2019}%
\special{pa 2892 2029}%
\special{pa 2924 2039}%
\special{pa 2955 2048}%
\special{pa 2987 2056}%
\special{pa 3019 2065}%
\special{pa 3051 2073}%
\special{pa 3117 2089}%
\special{pa 3149 2096}%
\special{pa 3215 2110}%
\special{pa 3249 2117}%
\special{pa 3282 2123}%
\special{pa 3315 2130}%
\special{pa 3349 2136}%
\special{pa 3382 2142}%
\special{pa 3450 2154}%
\special{pa 3480 2160}%
\special{fp}%
%
\special{pn 20}%
\special{pa 2224 1448}%
\special{pa 4776 1440}%
\special{fp}%
%
\special{pn 13}%
\special{pa 2008 1096}%
\special{pa 4792 1096}%
\special{dt 0.045}%
\special{pa 4792 1096}%
\special{pa 4792 1096}%
\special{dt 0.045}%
%
\special{pn 20}%
\special{pa 2000 1530}%
\special{pa 2350 1770}%
\special{fp}%
%
\special{pn 20}%
\special{pa 2000 1340}%
\special{pa 2290 1120}%
\special{fp}%
\put(19.1000,-17.3000){\makebox(0,0)[lb]{$N$}}%
\end{picture}}%

\caption{{1/2, 1/4 and 0 denote the $y$-coordinate, while 4
the $x$-coordinate.}}
\end{figure}
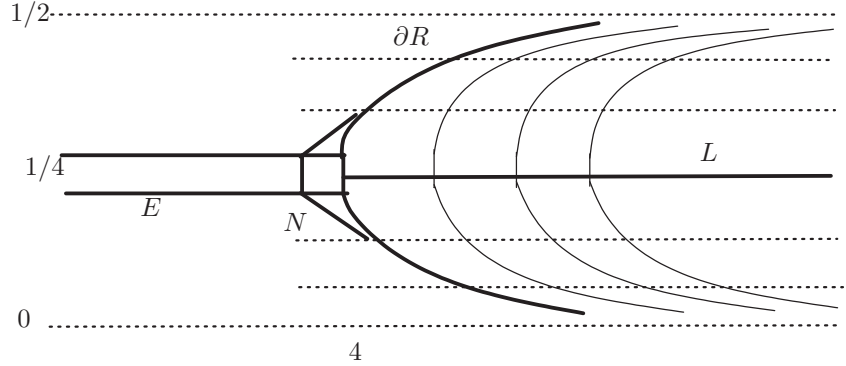

To restate Lemma \ref{lemma2}, we get:

\begin{corollary}
 \label{c}
The diffeomorphism $\phi$ is $1$-contracting along $T\FF^s$ on
$P'$, that is, if $v\in T_p\FF^s$, $p\in P'$, then
$\abs{D\phi(v)}\leq \abs{v}$. If furthermore $p\in N\cup R$, then
$\abs{D\phi(v)}= e^{-\lambda}\abs{v}$.
\end{corollary}
The strip $$E'=\{0\leq x\leq 4,\ \abs{y-\tfrac{1}{4}}<e^{-\lambda}\delta\}$$ is mapped to
the strip $$E=\{{2}\leq x\leq 4,\ \abs{y-\tfrac{1}{4}}<\delta\}$$ by a product map by virtue
of (\ref{12}).  Together with Lemma \ref{lemma}, we have:

\begin{lemma}\label{lemma3}
The foliation $\FF^u$ is vertical on the strip $E$.
\end{lemma}

We also have the following lemma by virtue of condition (F).
\begin{lemma}
 \label{lemma-added}
If $v$ is a vertical vector at a point on $E'$, then $\abs{D\phi(v)}=e^\lambda\abs{v}$.
\end{lemma}

\section{Expanding norm on $T\FF^u$}

Let 
$$\breve P=([0,\infty)\times[0,\tfrac{1}{2}])\cup\{x\geq 1-2y,\ \tfrac{1}{2}\leq y\leq 1\}.$$
See Figure 6. 
We shall define a metric $\Abs{\cdot}$ of $T\FF^u\vert_{\breve P}$
which is $e^\lambda$-expanding by $\phi$ in the sense that
$\Abs{D\phi(v)}\geq e^{\lambda}\Abs{v}$, 
$\forall v\in T\FF^u\vert_{\ \breve P}$.
The overall strategy is as follows. Suppose $\Abs{\cdot}_{\phi^{-1}(p)}$
is given. For any $v\in T\FF^u_{p}$, we shall define $\Abs{v}_p$ by
\begin{equation}
 \label{13}
\Abs{v}_p=\left\{ 
\begin{array}{ll}
e^\lambda\Abs{D\phi^{-1}(v)}_{\phi^{-1}(p)} & \mbox{if} \ \
 \tfrac{\abs{v}}{\abs{D\phi^{-1}(v)}}\leq e^\lambda\\
  \tfrac{\abs{v}}{\abs{D\phi^{-1}(v)}}\cdot \Abs{D\phi^{-1}(v)}_{\phi^{-1}(p)}   & 
\mbox{if} \ \ \tfrac{\abs{v}}{\abs{D\phi^{-1}(v)}}\geq e^\lambda,  
\end{array}\right.
\end{equation}
where $\abs{\cdot}$ denotes the Euclidean norm. 
\begin{figure}[h]
{\unitlength 0.1in%
\begin{picture}( 37.6800, 14.2600)(  6.0000,-18.1000)%
%
\special{pn 8}%
\special{pa 600 578}%
\special{pa 1736 1698}%
\special{fp}%
%
\special{pn 8}%
\special{pa 1728 1698}%
\special{pa 4368 1698}%
\special{fp}%
\special{pa 4368 1698}%
\special{pa 4360 1698}%
\special{fp}%
%
\special{pn 8}%
\special{pa 608 578}%
\special{pa 3472 578}%
\special{fp}%
%
\special{pn 20}%
\special{pa 1744 1674}%
\special{pa 1744 1138}%
\special{fp}%
\special{pa 1744 1138}%
\special{pa 1744 1138}%
\special{fp}%
%
\special{pn 20}%
\special{pa 608 562}%
\special{pa 616 562}%
\special{fp}%
\special{pa 616 578}%
\special{pa 1752 1154}%
\special{fp}%
\special{pa 1752 1154}%
\special{pa 1752 1154}%
\special{fp}%
%
\special{pn 8}%
\special{pa 1760 1154}%
\special{pa 4088 1154}%
\special{fp}%
%
\special{pn 20}%
\special{pa 3968 1706}%
\special{pa 3968 1154}%
\special{fp}%
\put(9.9200,-12.2600){\makebox(0,0)[lb]{$\Delta$}}%
\put(17.9200,-11.0600){\makebox(0,0)[lb]{$(0,1/2)$}}%
%
\put(17.9200,-11.0600){\makebox(0,0)[lb]{}}%
\put(6.3200,-5.1400){\makebox(0,0)[lb]{$(-1,1)$}}%
\put(42.8000,-12.1800){\makebox(0,0)[lb]{$\breve P$}}%
\put(13.2000,-8.5000){\makebox(0,0)[lb]{$X$}}%
\put(20.4000,-8.5000){\makebox(0,0)[lb]{$Y$}}%
%
\special{pn 20}%
\special{pa 2720 570}%
\special{pa 3460 830}%
\special{fp}%
\put(31.2000,-9.7000){\makebox(0,0)[lb]{$\phi(X)$}}%
\put(17.2000,-19.4000){\makebox(0,0)[lb]{$(0,0)$}}%
%
\special{pn 20}%
\special{pa 3450 840}%
\special{pa 3483 848}%
\special{pa 3516 857}%
\special{pa 3547 868}%
\special{pa 3576 881}%
\special{pa 3603 897}%
\special{pa 3627 917}%
\special{pa 3647 941}%
\special{pa 3664 968}%
\special{pa 3680 997}%
\special{pa 3695 1025}%
\special{pa 3713 1051}%
\special{pa 3734 1074}%
\special{pa 3757 1093}%
\special{pa 3784 1107}%
\special{pa 3813 1119}%
\special{pa 3845 1129}%
\special{pa 3877 1136}%
\special{pa 3911 1141}%
\special{pa 3946 1146}%
\special{pa 3980 1150}%
\special{fp}%
\end{picture}}%

\caption{}
\end{figure}
Let 
$$
X=(\{0\}\times[0,\tfrac{1}{2}])\cup\{x=1-2y,\ \tfrac{1}{2}\leq y\leq 1\}.$$
 We put the Euclidean norm on $X$ and
 apply the above strategy to get
a norm $\Abs{\cdot}$ on $\phi(X)$. Then interpolate in the region $Y$
bounded by $X$ and $\phi(X)$ the two norms monotonically
along the $\FF^s$-leaves. 
Apply the same strategy to $\phi(Y)$, and then to $\phi^2(Y)$ and so on.
Thus we obtain a norm on $\breve P\setminus R$.
But in fact, we can get a bit more. As is remarked in Lemma
\ref{lemma-added},
the map $\phi:E'\to E$ is already $e^\lambda$-expanding along
$T\FF^u$ with respect to
the Euclidean norm. Therefore
the norm we obtained on $E$ is nothing but the Euclidean norm.
Thus it extends continuously to $\partial R\cap E$, and one can
apply the same strategy {as in (\ref{13})} including this set. This way, we obtain a continuous 
norm on the closed set $\breve P\setminus {\rm Int}(R)$ which is
$e^\lambda$-expanding.
Next we shall extend the norm to $R$. Recall that the $\FF^u$-leaves in 
$R$ are the horizontal translates of $\partial R$. Define the norm on
each leaf simply as the translate of the norm $\Abs{\cdot}$
on $\partial R$. By the product structure of $R$
$$R\approx \partial R\times L$$ 
given by $\FF^u$ and
$\FF^s$, this norm on $R$ is also $e^\lambda$-expanding by $\phi$.

By Corollary
\ref{corollary}, the norm we obtained on the upper boundary of $\breve P$
is the image by $\tau$ of the norm on the lower boundary,
{as long as the interpolation in $Y$ is chosen to
be $\tau$-invariant on the horizontal boundaries}.
Therefore, by distributing the norm by the iterates of $\tau$,
 we get a continuous norm on $\cup_n\tau^n(\breve P)$ which is
$e^\lambda$-expanding by $\Phi$ and therefore by $F$.
Let 
$$\PP=\{x+y\geq0\}.$$
 Extend the norm $\Abs{\cdot}$ of
$T\FF^u$ from $\cup_n\tau^n(\breve P)$ to $\PP$ just setting {it} to be the
Euclidean norm on the difference set.
Summarizing the content of this setion, we get the following lemma.

\begin{lemma}
 \label{final-lemma}
There is a continous norm $\Abs{\cdot}$ on $T_{\PP}\FF^u$
with the following properties.

$\bullet$
$\Abs{v}\geq\abs{v}$
for any $v\in T_\PP\FF^u$,

$\bullet$ $\Abs{v}=\abs{v}$ for any $v\in T_\Delta\FF^u$,

$\bullet$
$\Abs{DF(v)}\geq e^\lambda\Abs{v}$ for any $v\in T_p\FF^u$,
 $p\in\cup_n\tau^n(\breve P)$.
\end{lemma}

\section{Final step}
We shall construct norms $\ABS{\cdot}$ along $T\FF^s$ and $T\FF^u$ on
$\PP$
for which $F$ is hyperbolic i.e, conditions (\ref{1}) and (\ref{2}) are
satisfied.
Recall that $T_p\FF^s$ is a horizontal line and $T_p\FF^u$ is a vertical
line for $p\in\Delta$ and that the differential $D\sigma$ of the
involution $\sigma$ maps $T_p\FF^s$ onto $T_p\FF^u$.
We shall construct $\ABS{\cdot}$ in such a way that 
\begin{equation}
 \label{symm}
\ABS{D\sigma(v)}=\ABS{v},  \ \ \ \ \ \forall v\in T_p\FF^s, \ \ \forall
p\in\Delta.
\end{equation}
Recall that $F^{-1}=\sigma F\sigma$, and we have $\sigma\FF^u=\FF^s$
and $\sigma\FF^s=\FF^u$.
After we have constructed the norms on $\mathbb P$, norms on
$\sigma(\mathbb P)$ will be given as the $D\sigma$-images.
That is,
$$\ABS{v}_p=\ABS{D\sigma(v)}_{\sigma(p)},\ \ \ \ \ \forall v\in T_p\FF^s\cup
T_p\FF^u, \ \ p\in\sigma(\PP).$$

Let $U=\{\abs{x+y}<1\}$, 
a { partial} fundamental domain of $F$.
We shall estimate the ratio $\ABS{DF(v)}/\ABS{v}$, 
$v\in T_p\FF^u\cup T_p\FF^s$, only when both $p$ and $F(p)$ are above $U$ or
below $U$. By the construction of $\ABS{\cdot}$ which follows, this ratio is
bounded when one of $p$ or $F(p)$ is contained in $U$.
To get the hyperbolicity, it is not a problem to skip one or two steps:
conditions (\ref{1}) and (\ref{2}) are asymptotic in nature.
Also hyperbolicity for the region below $U$ follows from the
hyperbolicity above $U$ by the symmetry.

Construction of $\ABS{\cdot}$ for $\PP\cap\{y<0\}$ is given in (I),
 and for $\PP\cap\{y>0\}$ in (II).
In (III), we shall show that the norms
constructed  yield a compete Riemannian metric.
Let $\epsilon$ be a positive number which is small compared with $\lambda$.

\bigskip
(I) {\sc Construction for $\PP_-=\PP\cap\{y<0\}$}.
 
For $v\in T_{(x,y)}\FF^s$, we let $\ABS{v}=e^{-\epsilon y}\abs{v}$.
By Corollary \ref{c}, $\Phi$ is 1-contracting along $T\FF^s$ with
respect
to the Euclidean metric and
$\tau$ is $e^{-\epsilon}$-contracting on $F^{-1}(\PP_-)\setminus U$ with
respect
to $\ABS{\cdot}$. Now it follows that 
$F$ is $e^{-\epsilon}$-contracting 
along $T\FF^s$ on this set. 
For $v\in T\FF^u$, define
$\ABS{v}=e^{-\epsilon y}\Abs{v}$. Then 
 $F$ is clearly $e^{\lambda-\epsilon}$-expanding along $T\FF^u$ on
 $F^{-1}(\PP)\setminus U$. 

Notice that the symmetry (\ref{symm}) is satisfied.
We do not estimate the contraction/expansion ratio on
$\PP\cap\{\abs{y}<1\}$ by the same reason as  we explained before. This is enough for
robust asymptotic estimates as in (\ref{1}).
Notice that $\ABS{v}\geq\abs{v}$ for $v\in T_p\FF^u\cup T_p\FF^s$,
$p\in\PP_-$.

\bigskip
(II) {\sc Construction for $\PP_+=\PP\cap\{y>0\}$}.

{Should we do the same construction as in (I) for the whole $\PP$,
an upward $\FF^u$-ray
 would have finite length, contrary to the completeness of the
metric. So we need a different construction
for $\PP_+$.}

 As for $T\FF^u$, we just put
$\ABS{\cdot}=\Abs{\cdot}$. Then $F$ is $e^\lambda$-expanding along
$T\FF^u$ on $\PP_+\setminus U$.

To define $\ABS{\cdot}$ on $T\FF^s$, {consider an arbitrary point} $p$ from the region
$C_n\subset\tau^n(\breve P)$
in Figure 7.
The point $p$ lies on a horizontal line segment which starts at a point 
$p_0\in\tau^n(X)$. 
Let $\ell$ be the distance between $p$ and $p_0$. Define 
$$\ABS{v}=e^{-\epsilon\ell}\abs{v}, \ \ \ \ v\in T_p\FF^s.$$
Next for a point $q$ from the region $D_n$, let $q_0\in B_n$ be the
point on the horizontal line passing through $q$. Define $\ABS{\cdot}$
on $T_q\FF^s$ to be equal to that on $T_{q_0}\FF^s$.
Here we make a natural identification of the horizontal line field:
$T\FF^s_q\approx T\FF^s_{q_0}$.
 Finally
on the subset $\PP_+\setminus\cup_{n\geq0}\tau^n(\breve P)$ (consisting
of small triangles), put $\ABS{\cdot}=\abs{\cdot}$.
Again the symmetry (\ref{symm}) is satisfied.

\begin{figure}[h]
{\unitlength 0.1in%
\begin{picture}( 36.8500, 12.9200)(  3.5700,-14.1000)%
%
\special{pn 8}%
\special{pa 357 311}%
\special{pa 4036 311}%
\special{fp}%
%
\special{pn 4}%
\special{pa 1352 1319}%
\special{pa 1352 1319}%
\special{fp}%
%
\special{pn 8}%
\special{pa 886 821}%
\special{pa 886 821}%
\special{dt 0.045}%
\special{pa 873 828}%
\special{pa 4042 821}%
\special{dt 0.045}%
\special{pa 4042 821}%
\special{pa 4030 802}%
\special{dt 0.045}%
%
\special{pn 8}%
\special{pa 1365 1313}%
\special{pa 1365 1313}%
\special{fp}%
%
\special{pn 8}%
\special{pa 1365 1332}%
\special{pa 3815 1332}%
\special{fp}%
%
\special{pn 20}%
\special{pa 1377 1344}%
\special{pa 1377 834}%
\special{fp}%
\special{pa 1377 834}%
\special{pa 1377 834}%
\special{fp}%
%
\special{pn 20}%
\special{pa 1377 834}%
\special{pa 376 317}%
\special{fp}%
\special{pa 376 317}%
\special{pa 376 317}%
\special{fp}%
%
\special{pn 20}%
\special{pa 2990 828}%
\special{pa 4023 324}%
\special{fp}%
%
\special{pn 20}%
\special{pa 3003 828}%
\special{pa 3003 979}%
\special{fp}%
%
\special{pn 20}%
\special{pa 3003 1168}%
\special{pa 3003 1332}%
\special{fp}%
%
\special{pn 20}%
\special{pa 3003 985}%
\special{pa 2959 996}%
\special{pa 2917 1007}%
\special{pa 2878 1019}%
\special{pa 2844 1030}%
\special{pa 2816 1041}%
\special{pa 2798 1053}%
\special{pa 2789 1065}%
\special{pa 2792 1077}%
\special{pa 2806 1090}%
\special{pa 2829 1102}%
\special{pa 2860 1115}%
\special{pa 2896 1128}%
\special{pa 2937 1141}%
\special{pa 2980 1154}%
\special{pa 3009 1162}%
\special{fp}%
%
\special{pn 8}%
\special{pa 3003 991}%
\special{pa 3034 984}%
\special{pa 3066 977}%
\special{pa 3159 956}%
\special{pa 3191 949}%
\special{pa 3222 943}%
\special{pa 3253 936}%
\special{pa 3285 930}%
\special{pa 3316 923}%
\special{pa 3348 917}%
\special{pa 3379 911}%
\special{pa 3411 906}%
\special{pa 3442 900}%
\special{pa 3474 895}%
\special{pa 3505 890}%
\special{pa 3537 885}%
\special{pa 3568 881}%
\special{pa 3664 869}%
\special{pa 3695 866}%
\special{pa 3791 857}%
\special{pa 3823 855}%
\special{pa 3855 852}%
\special{pa 3919 848}%
\special{pa 3929 847}%
\special{fp}%
%
\special{pn 8}%
\special{pa 3003 1168}%
\special{pa 3034 1176}%
\special{pa 3066 1183}%
\special{pa 3097 1191}%
\special{pa 3128 1198}%
\special{pa 3160 1205}%
\special{pa 3191 1213}%
\special{pa 3222 1220}%
\special{pa 3254 1226}%
\special{pa 3316 1240}%
\special{pa 3348 1246}%
\special{pa 3379 1252}%
\special{pa 3411 1257}%
\special{pa 3442 1263}%
\special{pa 3474 1268}%
\special{pa 3505 1272}%
\special{pa 3569 1280}%
\special{pa 3600 1283}%
\special{pa 3632 1286}%
\special{pa 3728 1292}%
\special{pa 3760 1293}%
\special{pa 3792 1293}%
\special{pa 3824 1294}%
\special{pa 3916 1294}%
\special{fp}%
\put(3.8800,-2.4800){\makebox(0,0)[lb]{$(-n-1,n+1)$}}%
\put(38.7200,-2.5400){\makebox(0,0)[lb]{$(n+1,n+1)$}}%
\put(32.1100,-6.2600){\makebox(0,0)[lb]{$B_n$}}%
\put(7.2200,-9.1600){\makebox(0,0)[lb]{$\Delta$}}%
\put(20.3900,-10.5400){\makebox(0,0)[lb]{$C_n$}}%
\put(37.1500,-7.1400){\makebox(0,0)[lb]{$D_n$}}%
%
\special{pn 8}%
\special{pa 804 538}%
\special{pa 1932 538}%
\special{dt 0.045}%
\special{pa 1932 538}%
\special{pa 1932 538}%
\special{dt 0.045}%
\put(12.5100,-4.7500){\makebox(0,0)[lb]{$\ell$}}%
%
\special{pn 4}%
\special{sh 1}%
\special{ar 1913 538 16 16 0  6.28318530717959E+0000}%
\special{sh 1}%
\special{ar 1919 538 16 16 0  6.28318530717959E+0000}%
\put(33.1500,-11.2000){\makebox(0,0)[lb]{$\tau^n(R)$}}%
\put(28.5000,-15.4000){\makebox(0,0)[lb]{$(n,n)$}}%
\put(17.9000,-4.6000){\makebox(0,0)[lb]{$e^{-\epsilon\ell}\vert\cdot\vert$}}%
\put(18.5000,-6.8000){\makebox(0,0)[lb]{$p$}}%
\put(7.7000,-4.7000){\makebox(0,0)[lb]{$p_0$}}%
\put(11.2000,-15.1000){\makebox(0,0)[lb]{$(-n,n)$}}%
\put(14.5000,-11.9000){\makebox(0,0)[lb]{$\tau^n(X)$}}%
%
\special{pn 4}%
\special{pa 380 310}%
\special{pa 1380 1330}%
\special{fp}%
\end{picture}}%

\caption{$B_n$ is a curve composed of line segments and a boundary curve
 of $\tau^n(R)$. When $n=0,1,2$, $B_n$ consists of just two line segments.
$D_n$ is the region right to $B_n$, Thus $\tau^n(R)\subset D_n$}
\end{figure}
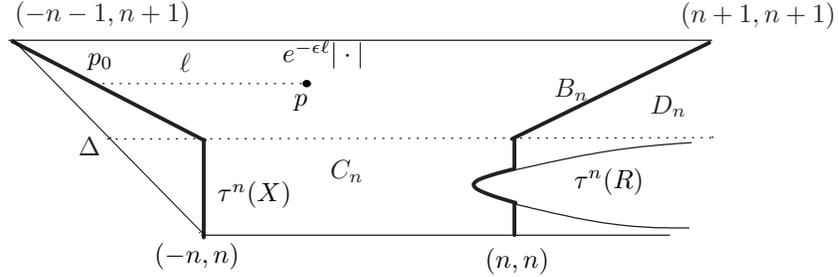

Now let us show that $F$ is $e^{-\alpha\epsilon}$-contracting along $T\FF^s$
on 
$$P_n'':=[-n,\infty)\times[n,n+1],$$
 where the constant $\alpha>0$ is from 
Lemma \ref{lemma2}. We assume $\alpha<1$.
First {consider the case where  $p$ lies on the upper half of $P''_n$:}
$p\in[-n,\infty)\times[n+\tfrac{1}{2},n+1]$.
There $\Phi$ is the $(2,0)$-translation and thus $F=\tau\circ\Phi$ is
the $(1,1)$-translation. {(The upper half of $P_0''$ is contained in the
region $C$ of Figure 3.) }
{If $p\in C_n$ (resp.\ $p\in D_n$), 
$f(p)\in C_{n+1}$ (resp.\ $f(p)\in D_{n+1}$). In both cases, we have
$$\ABS{DF(v)}=e^{-2\epsilon}\ABS{v},\ \ \ \forall v\in T_p\FF^s,$$
as is desired.}

{Next consider the case where $p$ lies in the lower half of 
$P_n''$ but not in the Reeb components $\tau^n(R)$. Notice that in this
part, the boundaries of $C_n$ are vertical,
and the norm $\ABS{\cdot}$ depends only on the $x$-coordinate.
 Thus
 in the computation of
contraction ratio, we only need to consider the function
$x\mapsto h_y(x)$: we do not have to care about the variation of
$y$-coordinate $y\mapsto g(y)$.

If $F(p)\not\in C_{n+1}$, then $F$ is $e^{-2\epsilon}$-contracting at $p$.
So consider the case $F(p)\in C_{n+1}$. 
If $p$ does not lie in $\tau^n(N)$,
then  $F$ is $e^{-\alpha\epsilon}$-contracting
 by virtue of Lemma \ref{lemma2}. 
If $p\in \tau^n(N)$, then $F$ is 
$e^{-\lambda}$-contracting, by virtue of
Corollary \ref{c} and the fact that $\Phi$ does not decrease the
$x$-coordinate in $N$. This holds true regardless of whether
$p\in C_{n}$ or not.}

For the Reeb component
$\tau^n(R)$, consider a horizontal ray $r$ contained in $\tau^n(R)$
with initial point on $\tau^n(\partial R)$.
The norm of $T_r\FF^s$ is determined by the $x$-coordinate of
the initial point.
Now the $x$-coordinate of
the initial point of $\Phi(r)$ is not less than the 
$x$-coordinate of
the initial point $q$ of $r$.
This shows that $F$ is $e^{-\lambda}$-contracting on $\tau^n(R)$
by Corollary \ref{c}.

It is clear that $\ABS{v}\geq\abs{v}$ for $v\in T_p\FF^u$, $p\in\PP_+$.
Our construction  $\ABS{\cdot}$ on $T_p\FF^s$ satisfies the following property, which
turns out to be useful in (III).

$\bullet$ {\em For any $n\in \N$,
there exists $c>0$ such that 
$\ABS{v}\geq c\abs{v}$ for $v\in T_p\FF^s$, $p\in\PP_+\cup\{y<n\}$.}

\bigskip
(III) We have defined the norm $\ABS{\cdot}$ on $\PP$.
As we said earlier, we define
the norm $\ABS{\cdot}$ on $\sigma(\PP)$ by transforming the former by
$D\sigma$. 
Define a Riemannian metric $m$ on $\R^2$ by using these norms and setting
that the two subspaces $T_p\FF^u$ and $T_p\FF^s$ be orthogonal.
We shall denote
by $\Abs{\cdot}_m$ the norm of $m$.
We have already shown that $F$ satisfies the hyperbolicity conditions
(\ref{1}) and (\ref{2}).
What is left is to show that $m$ is complete. 
Given arbitrarily large $R>0$, we shall show that the set $B(R)$ of points 
which are $R$-near to $(0,0)$ with respect to $m$ is bounded.
First given $n\in\Z$, consider the set 
$$
Q_n=\{n-\tfrac{1}{4}<y<n\}\cap\PP.$$
By Lemma \ref{lemma}, the foliation $\FF^u$ is vertical on $Q_n$ and
 $\Abs{v}_m\geq\abs{v}$ for any vertical
vector $v$ of $Q_n$.
This shows that any path in $\PP$ which crosses the strip $Q_n$ must have
$m$-length $\geq \tfrac{1}{4}$. The same is true for the strip
$\sigma(Q_n)$ in the region $\sigma(\PP)$. Thus the set $B(R)$ must be
contained in the region $Y$ bounded by $Q_{-n}\cup\sigma(Q_{-n})$
and $Q_n\cup\sigma(Q_n)$ for some $n>0$. But in $Y$, there is $c>0$ such that
$\Abs{v}_m\geq c\abs{v}$ for any tangent vector $v$ on $Y$.
In fact, if 
$$v=v_1+v_2, \ \ \ v_1\in T\FF^s, v_2\in T\FF^u, $$
we have
$\Abs{v}_m\geq\Abs{v_i}_m$ for each $i$
since $T\FF^u$ and $T\FF^s$ are orthogonal. On the other hand,
 there is $c>0$ depending on $Y$ such that $\Abs{v_i}_m\geq 2c\abs{v_i}$ for each $i$.
Now by the triangle inequality, there is $i$ such that
$\abs{v_i}\geq\abs{v}/2$. Then 
$$
\Abs{v}_m\geq\Abs{v_i}_m\geq 2c\abs{v_i}\geq c\abs{v},$$
as is desired. 
Now the set $B(R)$ must be contained in the Euclidean $R/c$--ball
centered at $(0,0)$. The proof of the completeness is now complete.

\bigskip
{\sc Final remark}. The diffeomorphism 
 $F$ is not topologically conjugate to a translation, since
the quotient space $\R^2/\langle F\rangle$ is not Hausdorff.
To show this, notice that any small piece of the $\FF^s$--leaf
passing through a point $p$ from the boundary of the Reeb component
$\partial R$ and
any small piece of the $\FF^u$--leaf passing through the point $\sigma(p)$
contain a common orbit.


\begin{thebibliography}{99}
\bibitem{GN} J. Groisman and Z. Nitecki, \em Foliations and
	conjugacy, II: The Mendes conjecture for time-one maps of flows,
	\rm arXiv:1812.04689.
\bibitem{M} P. Mendes, \em On Anosov diffeomorphisms on the plane, \em
	Proc. A.M.S. (1977) 231-235. 
\bibitem{W} W. White, \em An Anosov translation, \rm Dynamical Systems,
	    Proceedings of a Symposium held at the University of Bahia,
	    Salvador, Brasil, July 26-August 14, 1971  (M. M. Peixoto,
	    ed.) 1977, pp, 667-670.


\end{thebibliography}
\end{document}